\documentclass[hidelinks]{siamart251216}
\usepackage{amssymb,mathtools}
\usepackage{booktabs}
\usepackage{array}
\usepackage{enumitem}
\usepackage[protrusion=true,expansion=false]{microtype}
\usepackage{placeins}
\graphicspath{{figures/}}
\usepackage{bm}
\usepackage{textcomp}

\ifpdf
\hypersetup{
 pdftitle={Symbol-based residual certification of reduced and neural surrogates for tempered-fractional CGMY pricing},
 pdfauthor={D. Damircheli}}
\fi
\newsiamremark{remark}{Remark}
\newsiamthm{assumption}{Assumption}
\headers{Symbol-based residual certification for CGMY pricing}{D. Damircheli}
\title{Symbol-based residual certification of reduced and neural surrogates for tempered-fractional CGMY pricing}
\author{Davood Damircheli\thanks{Department of Mathematics, Louisiana State University, Baton Rouge, LA 70803, USA
(\email{ddamircheli@lsu.edu}).}}
\newcommand{\R}{\mathbb{R}}
\newcommand{\Pset}{\mathcal{P}}

\newcommand{\Vmu}{V_{\mu}}

\newcommand{\Mh}{M_h}
\newcommand{\Ah}{A_h}
\newcommand{\Sh}{S_h}
\newcommand{\norm}[1]{\left\lVert #1\right\rVert}

\newcommand{\dt}{\Delta\tau}
\newcommand{\etaSharp}{\eta_{\mathrm{sh}}}

\newcommand{\etaGuar}{\eta_{\mathrm{gu}}}
\newcommand{\ee}{\mathbf e}
\newcommand{\cc}{\mathbf c}
\newcommand{\rr}{\bm\rho}
\newcommand{\ff}{\mathbf f}
\newcommand{\vv}{\mathbf v}
\newcommand{\bb}{\mathbf b}
\newcommand{\zz}{\mathbf z}

\newcommand{\sres}{\mathbf s}
\newcommand{\dd}{\bm\delta}
\newcommand{\rf}{r}
\newcommand{\Hs}{H^{s_\mu}_{M,G}}
\newcommand{\Dm}{\mathcal D_G^-}
\newcommand{\Dp}{\mathcal D_M^+}
\newcommand{\Aop}{\mathcal A_\mu}
\DeclareMathOperator{\Real}{Re}
\newcommand{\Tn}[2]{T_{#1}(#2)}          %
\newcommand{\sym}[1]{#1}                 %
\newcommand{\smin}{s_{\min}}
\newcommand{\Thmu}{\Theta_\mu}

\begin{document}

\maketitle
\begin{abstract}
We certify reduced and neural approximations of European-put price trajectories under the
Carr--Geman--Madan--Yor model by residual bounds for the spline Galerkin discretisation of the
tempered-fractional pricing equation.  The Galerkin matrices are Toeplitz sections of an explicit
periodised Fourier symbol; summing the alias series to a stated tail makes the assembled system
the exact Galerkin system to $10^{-8}$.  A L\'evy--Khintchine argument gives the parameter-explicit
lower bound $\lambda_{\min}(S_I)\ge\min_\vartheta s_\mu(\vartheta)-\varepsilon_K\ge r\,m_{\min}$ for
the symmetric part, valid for $Y\in(0,2)\setminus\{1\}$, which replaces numerical spectral estimates
and positivity audits.  A banded comparison norm built from the mass and stiffness symbols yields a
rigorous certificate at $O(N)$ cost with explicit effectivity, and empirical interpolation of the
symbol family makes the operator affine in about 36 terms with a rigorous error transfer, so that a
reduced model is certified at a cost independent of $N$; a three-stage rule escalates to the
full-order bounds only when needed.  At the resolution that prices the whole parameter box to one
cent ($N=1024$), the full-order certificate costs $6\%$ of a direct solve, the certified
POD--Galerkin model of rank $80$ costs $3\%$ of it and settles $99.8\%$ of $10^4$ queries without
touching the full space.  No false accept occurs on any query with a reference solution (150 in
the parameter box, 264 outside it), and outside the box the certificate accepts 193 of the 210
queries on which the reduced model is still within one cent.  The neural
reduced models trained here fail the one-cent tolerance by a factor of at least $46$; the certificate
detects this at every query for a few milliseconds and explains it as a trade-off between
temporal smoothness and coefficient accuracy.  A residual certificate for the American put
in the same framework is proved as an outlook.
\end{abstract}

\section{Introduction}
\label{sec:introduction}
The CGMY tempered-stable process models asymmetric jumps and small-jump activity through four
parameters \cite{CarrGemanMadanYor2002}; European pricing under it leads to a nonlocal,
nonsymmetric parabolic PIDE whose fractional order varies with the model parameters
\cite{CarteadelCastilloNegrete2007}.  Repeated pricing, calibration and uncertainty
quantification motivate approximations of the parameter-to-price map: neural solution maps such as
DeepONet \cite{LuJinPangZhangKarniadakis2021} and Fourier neural operators
\cite{LiKovachkiEtAl2021,KovachkiEtAl2023}, nonlocal kernel networks
\cite{YouYuDEliaGaoSilling2022}, and neural reduced models that learn coordinates in a prescribed
trial space \cite{HesthavenUbbiali2018,BhattacharyaEtAl2021}.  An aggregate prediction error,
however, provides no error bound for an individual query outside well-represented regions of the
training set.  Query-wise residual bounds originate in certified reduced-basis methods for
parametric evolution problems \cite{GreplPatera2005,RozzaHuynhPatera2008,HaasdonkOhlberger2008},
where the coercivity constant is obtained by the successive constraint method
\cite{HuynhRozzaSenPatera2007} and non-affine operators by empirical interpolation
\cite{BarraultMadayNguyenPatera2004,GreplMadayNguyenPatera2007}; option-pricing applications
include \cite{ContLantosPironneau2011,BurkovskaHaasdonkSalomonWohlmuth2015} and parametric
Chebyshev interpolation of prices \cite{GassGlauMahlstedtMair2018}, and variationally correct
operator learning emphasises that the residual norm must reflect the stability of the underlying
problem \cite{QiuDahmenChen2026}.

This paper combines these ideas for CGMY pricing with a spline-valued surrogate: the deterministic
spline Galerkin discretisation \cite{DamircheliCGMY2026} defines the coefficient space, the training
trajectories, the residual and the correction operator (finite-difference, Galerkin and transform
methods \cite{ContVoltchkova2005,MatacheVonPetersdorffSchwab2004,HilberReichmannSchwabWinter2013,FangOosterlee2008}
supply the wider numerical context).  What is specific here is that the discretisation is
shift-invariant, so that the symbol calculus of Toeplitz sections
\cite{GrenanderSzego1958,BottcherSilbermann1999} replaces the successive constraint method by a
closed-form constant, the interpolation acts on a function of one variable, and the circulant
preconditioner of Chan and Strang \cite{ChanStrang1989} makes the exact Riesz solve a dozen
iterations at any resolution.  The contributions are:
\begin{enumerate}[leftmargin=2em]
\item exact assembly of the spline Galerkin matrices from the periodised
Fourier symbol with an explicit truncation bound (Lemmas~\ref{lem:periodisation}--\ref{lem:tail});
\item a rigorous, parameter-explicit spectral lower bound for the free symmetric part
(Theorem~\ref{thm:spectral}), valid for $Y\in(0,2)\setminus\{1\}$, replacing the
numerical spectral estimate and the positivity audit of the earlier version;
\item a banded comparison norm giving an $O(N)$ residual certificate with explicit effectivity,
and a two-stage rule with a circulant-preconditioned exact Riesz solve
(Theorem~\ref{thm:comparison}, Proposition~\ref{prop:twostage});
\item an affine decomposition of the symbol family with rigorous error transfer, hence an
$N$-independent certified reduced model (Lemma~\ref{lem:eim}, Theorem~\ref{thm:eim-cert});
\item matched evaluations at the resolution that prices the parameter box to one cent: the cost
of the certificate at full and reduced order, a certified reduced model an order of magnitude
faster than the direct solve, out-of-distribution detection, and a negative result for the neural
reduced models that the certificate itself establishes and explains.
\end{enumerate}
The discrete price is compared with the COS reference \cite{FangOosterlee2008} throughout, so that
the certified trajectory bound is a stated part of a one-cent price budget rather than a substitute
for it.  The principal computational finding is that the guarantee is cheap and the reduced model
earns it: at $N=1024$ the full-order certificate costs $6\%$ of the direct solve, the certified
reduced model of rank $80$ runs at $3\%$ of it with $99.8\%$ of $10^4$ queries settled at the
$N$-independent stage, while the neural reduced models trained here fail the tolerance by a
factor of at least $46$ and are rejected at every query for a few milliseconds.

The analysis distinguishes the assembled discrete system from the exact Galerkin discretisation
and the continuum price: residual certification controls the distance to the exact solution of the
first, the variational results describe the second under exact assembly, and Fourier comparisons
assess the third empirically; the implementation uses ordinary floating point, with its arithmetic
assumptions stated with the bracket.  Section~\ref{sec:model} defines the model and the spline
Galerkin discretisation via symbols, Section~\ref{sec:certificate} develops certification and
correction, Section~\ref{sec:affine} the affine structure and the candidate approximations,
Section~\ref{sec:numerics} evaluates them, and Section~\ref{sec:conclusion} interprets the
findings; the supplementary materials give the uniform-stability proof with the explicit constant,
the verification studies, supporting extensions and the American-put certificate.  Throughout,
$\psi$ denotes the payoff, $g$ its exterior lifting, $\mu=(C,G,M,Y,r)$ the spatial-operator
parameters, and $T$ the maturity.
Table~\ref{tab:notation} collects the symbols used throughout.
\begin{table}[ht]
\centering\footnotesize
\caption{Notation.}
\label{tab:notation}
\begin{tabular}{@{}l@{\ \ }>{\raggedright\arraybackslash}p{3.3cm}@{\ \ }l@{\ \ }>{\raggedright\arraybackslash}p{3.3cm}@{}}
\toprule
$\mu=(C,G,M,Y,r)$, $T$, $\mathcal K$ & model parameters, maturity, strike &
$N$, $h$, $p$ & cells, mesh width, spline order ($p=3$)\\
$N_h$, $N_T$, $\dt$ & free unknowns, time steps, step size &
$A_\mu(\xi)$ & full-line symbol of $\Aop$, \eqref{eq:symbol}\\
$\sym{a}_\mu$, $\sym{s}_\mu=\Real\sym{a}_\mu$ & periodised symbol, its real part &
$\sym{m},\sym{d}_1,\sym{g}_1$ & mass, drift, stiffness symbols\\
$K$, $\varepsilon_K$ & alias truncation, tail bound &
$\Tn{n}{f}$ & $n\times n$ Toeplitz section of $f$\\
$M_I,A_I,S_I$ & free mass, operator, symmetric part &
$\alpha_{\rm lb}$, $\smin$ & spectral lower bound, $\min_\vartheta\sym{s}_\mu$\\
$X_I$, $\beta_\mu$, $\Thmu$ & banded comparison norm, its factors &
$\bar X$, $\bar\Theta_\mu$ & cell comparison norm (reduced stage)\\
$\rr^{n+1/2}$, $\beta_0$ & CN residual, initial discrepancy &
$\eta^S,\eta^X,\eta^{\bar X}$ & exact, banded, reduced estimators\\
$Q$, $\varepsilon_Q$, $A_Q$ & EIM rank, error, affine operator &
$r$, $V$, $\varrho$ & POD rank, basis, retained rows\\
$\varepsilon$, $\varepsilon_N$ & tolerance; one cent at $N=1024$ &
$\eta_{\rm aug}$, $\lambda^{\bar X}_{\rm aug}$ & augmented bound, lower companion\\
\bottomrule
\end{tabular}
\end{table}
\section{CGMY model and constrained discretisation}
\label{sec:model}\label{sec:galerkin}
\subsection{CGMY backward operator}
\label{ssec:operator}

Let $S_t=e^{X_t}$ under a risk-neutral measure and let $V(x,\tau)$ be the
option value in the log-price variable $x=\log S$ and backward time
$\tau=T-t$.  Following the convention of the companion paper
\cite{DamircheliCGMY2026}, the CGMY parameters satisfy
\begin{equation}
  C>0,\qquad G>0,\qquad M>1,\qquad Y\in(0,2),\qquad Y\ne1,
  \label{eq:cgmy-parameters}
\end{equation}
where $M>1$ is the exponential-moment condition that makes the stock
martingale correction
\begin{equation}
  \omega_S(\mu)
  = C\Gamma(-Y)\bigl[(M-1)^Y-M^Y+(G+1)^Y-G^Y\bigr]
  \label{eq:martingale-correction}
\end{equation}
real valued.  Set
\begin{equation}
  \lambda_0(\mu) := \rf + C\Gamma(-Y)\bigl(M^Y+G^Y\bigr).
  \label{eq:lambda0}
\end{equation}
For $Y\in(1,2)$ one has $\Gamma(-Y)>0$ and hence $\lambda_0>\rf$.  The
tempered left- and right-sided Riemann--Liouville operators are
\begin{equation}
  \Dm u=e^{-Gx}\,{}_{-\infty}D_x^Y\bigl(e^{Gx}u\bigr),
  \qquad
  \Dp u=e^{Mx}\,{}_{x}D_\infty^Y\bigl(e^{-Mx}u\bigr),
  \label{eq:tempered-ops}
\end{equation}
and the single backward-time operator used for both the theory and the
discretisation is
\begin{equation}
  \Aop u
  =-(\rf-\omega_S)u_x
   -C\Gamma(-Y)\Dp u
   -C\Gamma(-Y)\Dm u
   +\lambda_0u.
  \label{eq:cgmy-operator}
\end{equation}
The European pricing problem is
\begin{equation}
  \partial_\tau V+\Aop V=0,
  \qquad
  V(\cdot,0)=\psi.
  \label{eq:pricing-pide}
\end{equation}
On the full line the Fourier symbols are
$\widehat{\Dm u}(\xi)=(G+i\xi)^Y\widehat u(\xi)$ and
$\widehat{\Dp u}(\xi)=(M-i\xi)^Y\widehat u(\xi)$ on the principal branch, so
that
\begin{equation}
  A_\mu(\xi)
  =-(\rf-\omega_S)i\xi
   -C\Gamma(-Y)\bigl[(M-i\xi)^Y+(G+i\xi)^Y\bigr]
   +\lambda_0 .
  \label{eq:symbol}
\end{equation}

\subsection{Parameter family and exterior data}
\label{ssec:family}\label{ssec:finite-domain}
The coefficient parameters lie in a compact interior box
\begin{equation}
\begin{aligned}
 \Pset=\bigl\{\mu=(C,G,M,Y,\rf):\;&
 C\in[C_{\min},C_{\max}],\ G\in[G_{\min},G_{\max}],\ M\in[M_{\min},M_{\max}],\\
 &Y\in[Y_{\min},Y_{\max}],\ \rf\in[r_{\min},r_{\max}]\bigr\}
\end{aligned}
 \label{eq:parameter-set}
\end{equation}
with $C_{\min}>0$, $G_{\min}>0$, $M_{\min}>1$, $[Y_{\min},Y_{\max}]\Subset(1,2)$ and $r_{\min}>0$;
the first, fourth and fifth restrictions give uniform coercivity and $M_{\min}>1$ selects the
interior exponential-moment regime.  Maturity $T$ is an additional surrogate input and $\psi$ the
payoff.  The numerical study uses the box
\begin{equation}
\begin{gathered}
 C\in[0.05,1],\quad G\in[1,5],\quad M\in[2,9],\quad Y\in[1.1,1.8],\quad \rf\in[0.01,0.1],\\
 T\in[0.1,1]\ \text{years},\quad \mathcal K=100,
\end{gathered}
 \label{eq:numerical-box}
\end{equation}
and a common truncated domain $\Omega=(x_{\min},x_{\max})$ throughout the family.

Nonlocal exterior values enter through $V=g+U$, with $U$ extended by zero
outside $\Omega$. For a put, a smooth lifting is
$g(x,\tau)=(\mathcal Ke^{-r\tau}-e^x)\zeta(x)$, where $\zeta$ equals one to the left
of the domain and zero to its right. With $s_\mu=Y/2$, define
\begin{equation}
  \Vmu=\Hs(\Omega)
  =\bigl\{v\in L^2(\Omega):\widetilde v\in H^{s_\mu}(\R)\bigr\},
  \qquad
  \norm{v}_{\Vmu}:=\norm{\widetilde v}_{H^{s_\mu}_{M,G}(\R)},
  \label{eq:energy-space}
\end{equation}
where $\widetilde v$ is extension by zero and the full-line norm is defined by
\begin{equation}
  \norm{u}_{H^{q}_{M,G}(\R)}^2
  :=\frac1{2\pi}\int_\R
  \Bigl(1+|M-i\xi|^{2q}+|G+i\xi|^{2q}\Bigr)|\widehat u(\xi)|^2\,d\xi .
  \label{eq:HMG-norm}
\end{equation}
The finite-domain bilinear form is
$a_\mu(u,w):=\langle\Aop\widetilde u,\widetilde w\rangle$, and the truncated
weak problem is: find $U\in L^2(0,T;\Vmu)$ with
$\partial_\tau U\in L^2(0,T;\Vmu')$, $U(\cdot,0)=\psi-g(\cdot,0)$, and
\begin{equation}
  (\partial_\tau U,w)_{L^2(\Omega)}+a_\mu(U,w)=\ell^g(w;\tau)
  \qquad\forall\,w\in\Vmu,
  \label{eq:weak-problem}
\end{equation}
where the lifting functional is
\begin{equation}
  \ell^{g}(w;\tau)
  :=-\int_\Omega\partial_\tau g\,w\,dx-\int_\Omega(\Aop g)\,w\,dx .
  \label{eq:lifting-functional}
\end{equation}
All exterior data enter only through $\Aop g$ evaluated on $\Omega$.

The put payoff is $\psi(x)=(\mathcal K-e^x)^+$ with strike $\mathcal K$ (the letter $K$ is reserved
for the alias truncation of Section~\ref{ssec:assembly}). The main study fixes $\mathcal K=100$ and
learns variation in $(\mu,T)$; this is a parameter-to-trajectory map within
one payoff family.

\subsection{Uniform variational stability}
\label{ssec:uniform}
\begin{theorem}[Parameter-uniform continuity and coercivity]
\label{thm:uniform-coercivity}
For the family \eqref{eq:parameter-set}, there are constants
$C_{\rm cont}<\infty$ and $c_1>0$, independent of $\mu$, such that
\begin{equation}
 |a_\mu(v,w)|\le C_{\rm cont}\|v\|_{V_\mu}\|w\|_{V_\mu},
 \qquad a_\mu(v,v)\ge c_1\|v\|_{V_\mu}^2
 \quad(v,w\in V_\mu).
 \label{eq:uniform-coercivity}
\end{equation}
\end{theorem}
\begin{proof}
Explicit choices and the proof appear in Section~\ref{app:proof-uniform} of the supplementary materials.
\end{proof}
The constants establish uniform stability despite the varying energy order.
They are conservative; the discrete bounds of Section~\ref{ssec:spectral} are the ones used
in the certificate.  The continuum statement keeps the hypothesis $Y>1$: for $Y<1$ the drift
term is not controlled by the $H^{Y/2}_{M,G}$ norm.  The discrete results below hold for
$Y\in(0,2)\setminus\{1\}$.
\begin{remark}[Explicit coercivity constant]\label{rem:coercivity}
Lemma~\ref{lem:LK} also sharpens Theorem~\ref{thm:uniform-coercivity}.  Using $1-\cos u\ge2u^2/\pi^2$
on $|u|\le\pi$,
\[
 I_a(\xi)\ \ge\ \frac{2}{\pi^2}\gamma_Y\,a^{Y-2}\xi^2\ \ (|\xi|\le a),\qquad
 I_a(\xi)\ \ge\ \frac{2e^{-\pi}}{(2-Y)\pi^Y}\,|\xi|^Y\ \ (|\xi|\ge a),
\]
with $\gamma_Y:=\int_0^\pi u^{1-Y}e^{-u}du$, which yield the explicit constant
$c_1=5.1\times10^{-5}$ of \eqref{eq:c1-explicit} on the box \eqref{eq:numerical-box} with a
half-page proof (Section~\ref{app:proof-uniform} of the supplementary materials), against
$1.6\times10^{-8}$ for the earlier high-frequency argument; the optimal constant, attained at
$\xi=0$, is $r_{\min}/(1+M_{\max}^{Y_{\max}}+G_{\max}^{Y_{\max}})=1.40\times10^{-4}$.
\end{remark}

\subsection{Spline space and assembly}
\label{ssec:spline-space}
Let $N_p$ be the cardinal spline of order $p$ (degree $p-1$),
$h=(x_{\max}-x_{\min})/N$, and
$\phi_{h,k}(x)=h^{-1/2}N_p((x-x_{\min})/h-k)$. The homogeneous space is
\begin{equation}
 V_h^0=\operatorname{span}\{\phi_{h,k}:k=0,\ldots,N-p\},
 \qquad N_h=N-p+1.
 \label{eq:Vh}
\end{equation}
For $p\ge2$, these interior-supported splines have zero extensions in
$H^1(\R)$ and hence in $V_\mu$ for $Y/2<1$. Exact variational assembly gives
\begin{equation}
 (\Mh)_{kk'}=(\phi_{h,k'},\phi_{h,k})_{L^2(\Omega)},
 \qquad
 (\Ah(\mu))_{kk'}=a_\mu(\phi_{h,k'},\phi_{h,k}),
 \label{eq:matrices}
\end{equation}
so that
\begin{equation}
 \Ah(\mu)
 =-(\rf-\omega_S)D_h
  -C\Gamma(-Y)T_M^+
  -C\Gamma(-Y)T_G^-
  +\lambda_0\Mh,
 \label{eq:Ah}
\end{equation}
where $D_h$ is the first-derivative matrix and the fractional blocks are
Toeplitz with the $h^{-Y}$ scaling
\begin{equation}
 (T_G^-)_{kk'}=h^{-Y}\Lambda_G^-(k-k';Gh),
 \qquad
 (T_M^+)_{kk'}=h^{-Y}\Lambda_M^+(k-k';Mh),
 \label{eq:toeplitz-blocks}
\end{equation}
the tempered connection coefficients $\Lambda_G^-,\Lambda_M^+$ being the absolutely convergent
Fourier-symbol integrals made explicit in Section~\ref{ssec:assembly}.  The matrices $\Mh$ and $D_h$
are banded; the Toeplitz blocks are applied by circulant embedding in
$O(N_h\log N_h)$ operations with no compression error.  The semidiscrete
system is
\begin{equation}
  \Mh\dot\cc(\tau)+\Ah(\mu)\cc(\tau)=\ff(\tau;\mu,\psi),
  \qquad
  f_k(\tau)=\ell^{g}(\phi_{h,k};\tau).
  \label{eq:semidiscrete-system}
\end{equation}

\subsection{Exact assembly from the periodised symbol}
\label{ssec:assembly}
Every matrix in \eqref{eq:Ah} is the Galerkin matrix of a Fourier multiplier in a shift-invariant
basis, and is therefore a Toeplitz matrix whose entries are the Fourier coefficients of a
$2\pi$-periodic function, the \emph{periodised symbol}.  This observation replaces quadrature of
the connection coefficients by an alias sum with an explicit tail.
Write $\vartheta_k:=\vartheta+2\pi k$ and
$\widehat N_p(\vartheta)=\bigl(\tfrac{\sin(\vartheta/2)}{\vartheta/2}\bigr)^p$
(up to a unimodular factor), so that $|\widehat N_p(\vartheta)|^2\le(2/|\vartheta|)^{2p}$.
\begin{lemma}[Periodised symbol]\label{lem:periodisation}
Let $p\ge2$ and $Y\in(0,2)\setminus\{1\}$.  The exact Galerkin matrix of
$\Aop$ in the basis $\{\phi_{h,k}\}$ is the Toeplitz matrix
$(\Ah)_{kk'}=t_{k-k'}$ with
\begin{equation}
\begin{gathered}
 t_m=\frac1{2\pi}\int_{-\pi}^{\pi}\sym{a}_\mu(\vartheta)e^{i\vartheta m}\,d\vartheta,
 \\
 \sym{a}_\mu(\vartheta):=\sum_{k\in\mathbb Z}\bigl|\widehat N_p(\vartheta_k)\bigr|^2A_\mu(\vartheta_k/h),
\end{gathered}
 \label{eq:periodised}
\end{equation}
where the series converges absolutely and uniformly on $[-\pi,\pi]$, and
$\sym{a}_\mu$ is real-analytic and $2\pi$-periodic.  Likewise the mass matrix,
the first-derivative matrix and the $H^1$ stiffness matrix $h^2R^{(1)}_h$
have symbols
$\sym{m}(\vartheta)=\sum_k|\widehat N_p(\vartheta_k)|^2$,
$\tfrac{i}{h}\sym{d}_1(\vartheta)$ with $\sym{d}_1(\vartheta)=\sum_k|\widehat N_p(\vartheta_k)|^2\vartheta_k$,
and $\sym{g}_1(\vartheta)=\sum_k|\widehat N_p(\vartheta_k)|^2\vartheta_k^2$.
For $p=3$,
$\sym{m}(\vartheta)=\tfrac{11}{20}+\tfrac{13}{30}\cos\vartheta+\tfrac1{60}\cos2\vartheta$,
$\sym{d}_1(\vartheta)=\tfrac56\sin\vartheta+\tfrac1{12}\sin2\vartheta$,
$\sym{g}_1(\vartheta)=1-\tfrac23\cos\vartheta-\tfrac13\cos2\vartheta$,
and $m_{\min}:=\min_\vartheta \sym{m}=\sym{m}(\pi)=2/15$.
\end{lemma}
\begin{proof}
With $\widehat\phi_{h,k}(\xi)=h^{1/2}\widehat N_p(h\xi)e^{-i\xi(x_{\min}+kh)}$ and
Plancherel,
\[
a_\mu(\phi_{h,k'},\phi_{h,k})=\frac1{2\pi}\int_\R A_\mu(\xi)\,h|\widehat N_p(h\xi)|^2e^{i\xi h(k-k')}d\xi;
\]
substitute $\vartheta=h\xi$ and split $\R=\bigcup_k[-\pi,\pi)+2\pi k$.  Since
$|A_\mu(\xi)|\lesssim 1+|\xi|^{\max(Y,1)}$ and $|\widehat N_p(\vartheta)|^2\lesssim|\vartheta|^{-2p}$,
the terms are bounded by $|\vartheta_k|^{\max(Y,1)-2p}$, summable because
$2p-\max(Y,1)>1$.  Analyticity follows from that of each term and uniform
convergence of derivatives.  The mass, drift and stiffness symbols are
trigonometric polynomials of degree $p-1$ whose coefficients are the exact
entries $\int N_pN_p(\cdot-m)$, $\int N_p'N_p(\cdot-m)$, $\int N_p'N_p'(\cdot-m)$
(for $p=3$: $\tfrac1{120}(1,26,66,26,1)$, $(\tfrac1{24},\tfrac5{12},0,-\tfrac5{12},-\tfrac1{24})$
and $(-\tfrac16,-\tfrac13,1,-\tfrac13,-\tfrac16)$).
\end{proof}

\begin{lemma}[Truncation tail]\label{lem:tail}
Let $\sym{a}_{\mu,K}$ be the series \eqref{eq:periodised} truncated to $|k|\le K$,
and for $q\in[0,2p-1)$ put
$\tau_K(q):=\dfrac{2^{2p+1}\pi^{q-2p}(2K-1)^{\,q-2p+1}}{2(2p-q-1)}$.
Then for all $|\vartheta|\le\pi$,
\begin{equation}
\begin{aligned}
 |\sym{a}_\mu(\vartheta)-\sym{a}_{\mu,K}(\vartheta)|
 \le \varepsilon_K(\mu):={}&\frac{|\rf-\omega_S|}{h}\,\tau_K(1)+\lambda_0\,\tau_K(0)\\
 &+2C|\Gamma(-Y)|h^{-Y}\Bigl(1+\frac{\max(G,M)h}{\pi(2K+1)}\Bigr)^{Y}\tau_K(Y).
\end{aligned}
 \label{eq:tail}
\end{equation}
If only the fractional part of \eqref{eq:periodised} is truncated and the
trigonometric polynomials $\sym{m},\sym{d}_1$ are used exactly, the first and
third terms are absent.
\end{lemma}
\begin{proof}
For $|k|\ge K+1$ and $|\vartheta|\le\pi$, $|\vartheta_k|\ge\pi(2|k|-1)\ge\pi(2K+1)$.
Hence $|\widehat N_p(\vartheta_k)|^2|\vartheta_k|^q\le2^{2p}(\pi(2|k|-1))^{q-2p}$ and
$|\delta+i\vartheta_k|^Y\le|\vartheta_k|^Y(1+\delta/(\pi(2K+1)))^Y$.  Sum over
$|k|\ge K+1$ using $\sum_{k\ge K+1}(2k-1)^{-s}\le\int_K^\infty(2t-1)^{-s}dt$.
\end{proof}

The series \eqref{eq:periodised} is evaluated on the grid $\vartheta_j=-\pi+2\pi j/L$ with
$L\ge 2N_h$ a multiple of $N_h$; the $2N_h-1$ Toeplitz entries then follow from one inverse FFT,
$t_m=(-1)^m\,[\mathcal F_L^{-1}\sym{a}_{\mu,K}]_{m\bmod L}$, with an aliasing error
$\sum_{j\ne0}|t_{m+jL}|$ that decays like $e^{-\min(G,M)h\,(L-N_h)}$ because the fractional symbol is
analytic in the strip $|\Im\vartheta|<\min(G,M)h$; with $L=4N_h$ the exponent is at least $46$ on the
box.  The cost is $O(LK)$ scalar operations per query.
All experiments use $K=100$ with the exact trigonometric polynomials for the mass, drift and
stiffness parts; then $\varepsilon_K(\mu)\le5\times10^{-9}\,\norm{\Ah}_2$ on the box of
\eqref{eq:numerical-box} (Table~\ref{tab:assembly-v9}), and $\varepsilon_K$ is recorded with
every assembled operator because it enters the certificate.  The principal-period rule of the
earlier version of this work, which samples only the $k=0$ term, is the case $K=0$; its relative
operator-norm error is $0.31$--$0.41$ at every resolution because at $\vartheta=\pm\pi$ the $k=0$
and $k=\mp1$ terms of \eqref{eq:periodised} have equal modulus.  With $K\ge100$ the assembled
matrices are the exact Galerkin matrices to $10^{-8}$, and every exact-assembly statement applies
to them with the explicit perturbation \eqref{eq:tail}; \eqref{eq:periodised} is checked
independently against entries computed at 40 digits (Section~\ref{ssec:e0-assembly}).
The algebraic facts about Toeplitz sections that the certificate uses are classical
\cite{GrenanderSzego1958,BottcherSilbermann1999}; we record them because the constants matter.
\begin{lemma}[Symbol inequalities]\label{lem:toeplitz}
For bounded $2\pi$-periodic $f,g$ with $g>0$ real and every $v\in\R^n$,
\[
 v^\top\Tn{n}{f}v=\frac1{2\pi}\int_{-\pi}^{\pi}f(\vartheta)\Bigl|\sum_kv_ke^{ik\vartheta}\Bigr|^2d\vartheta,
 \qquad
 \norm{\Tn{n}{f}}_2\le\norm f_\infty,
\]
and
\[
 v^\top\Tn{n}{f}v\ \ge\ \Bigl(\inf_\vartheta\frac{f}{g}\Bigr)\,v^\top\Tn{n}{g}v .
\]
In particular $\lambda(\Tn{n}{f})\subset[\inf f,\sup f]$ for real $f$.
\end{lemma}
\begin{proof}
Expand the quadratic form; the norm bound is the compression of
multiplication by $f$ on $L^2(\mathbb T)$; the last inequality is the first
identity applied to $f-(\inf f/g)\,g\ge0$.
\end{proof}

\subsection{Prescribed coefficients and the free system}
\label{ssec:constrained-system}
Write $\bar M,\bar A$ for the raw $(N+p-1)$-dimensional assembled matrices, let $E$ inject the
$N_h=N-p+1$ free coefficients into the raw vector, and let $b^n$ vanish on the free entries and
carry the prescribed boundary coefficients: the left ones equal $\mathcal K e^{-r\tau_n}\sqrt h$,
the right ones zero.  Every admissible trajectory has the form $c^n=Ez^n+b^n$, and
\[
 M_I=E^\top\bar M E,\quad A_I=E^\top\bar A E,\quad
 S_I=\tfrac12(A_I+A_I^\top)=E^\top\bar S E,\quad \bar S=\tfrac12(\bar A+\bar A^\top).
\]
These are coefficient constraints, not extra equations.  Restricting the Crank--Nicolson (CN)
equations to the free rows gives
\begin{equation}
 M_I\frac{z^{n+1}-z^n}{\Delta\tau}
 +A_I\frac{z^{n+1}+z^n}{2}=f_I^{n+1/2},
 \label{eq:constrained-cn}
\end{equation}
with the effective free forcing
$f_I^{n+1/2}=E^\top[\bar f^{n+1/2}-\bar M\,(b^{n+1}-b^n)/\Delta\tau-\bar A\,(b^{n+1}+b^n)/2]$,
and for a candidate with the identical $b^n$ the free rows of the raw residual,
\begin{equation}
 \rho_I^{n+1/2}=E^\top\!\left[\bar M\frac{c_\theta^{n+1}-c_\theta^n}{\Delta\tau}
 +\bar A\frac{c_\theta^{n+1}+c_\theta^n}{2}-\bar f^{n+1/2}\right],
 \label{eq:constrained-residual}
\end{equation}
are exactly the lifted free-space residual; boundary discrepancies must first be removed by
projection, and the error satisfies $e^n=Ee_I^n$, $(e^n)^\top\bar M e^n=(e_I^n)^\top M_I e_I^n$.
Because the free indices are contiguous, the free matrices are Toeplitz sections of the
symbols of Lemma~\ref{lem:periodisation}:
\begin{equation}
 A_I=\Tn{N_h}{\sym{a}_\mu},\quad M_I=\Tn{N_h}{\sym{m}},\quad
 S_I=\Tn{N_h}{\sym{s}_\mu},\quad \sym{s}_\mu:=\Real\sym{a}_\mu,\quad
 h^2R^{(1)}_I=\Tn{N_h}{\sym{g}_1},
 \label{eq:sections}
\end{equation}
the drift symbol being odd and therefore absent from the symmetric part; $h^2R^{(1)}_I$ is banded with
bandwidth $p-1$.  The raw matrices are the sections $\Tn{N+p-1}{\sym{m}}$ and $\Tn{N+p-1}{\sym{a}_\mu}$ of the
full cardinal splines, so the prescribed boundary functions are not truncated to $\Omega$ in the
lifting term; the initial projection $c^0_h$ uses the $\Omega$-quadrature mass matrix so that
constants are reproduced exactly up to the boundary.  At $N=64$ and $p=3$ (quadratic splines) there
are 66 stored coefficients and 62 free unknowns.  Certification uses $M_I$ itself and restricts all
operator and preconditioner actions to the free space; positivity of $S_I$ is a theorem
(Section~\ref{ssec:spectral}).

\subsection{Time discretisation and error target}
\label{ssec:reference-time}\label{ssec:deterministic-error}\label{ssec:lifting-cancels}
On the uniform grid $\tau_n=n\dt$, $\dt=T/N_T$, write $\overline v^{n+1/2}=(v^{n+1}+v^n)/2$.  The
reference uses pure CN without startup smoothing.  In the algebraic analysis
the tuple $(\Mh,\Ah,\Sh,\cc,\ff)$ denotes $(M_I,A_I,S_I,\zz,f_I)$, and in exact-assembly statements
their variational counterparts; the reference solves
\begin{equation}
 \Mh\frac{\cc_h^{n+1}-\cc_h^n}{\dt}
 +\Ah(\mu)\overline\cc_h^{\,n+1/2}
 =\overline\ff^{\,n+1/2}.
 \label{eq:cn-reference}
\end{equation}
Reference and candidate share the lifting $b^n$, so
$\|V_h^n-V_{\theta,h}^n\|_{L^2(\Omega)}=\|\ee^n\|_{\Mh}$ for the coefficient error $\ee^n$, and the
continuum error decomposes as
\begin{equation}
 V(\tau_n)-V_{\theta,h}^n
 =[V(\tau_n)-V_h^n]+[V_h^n-V_{\theta,h}^n].
 \label{eq:total-error-split}
\end{equation}
The certificate controls the second term; the first admits, under the regularity and buffer
hypotheses of \cite{DamircheliCGMY2026}, the deterministic estimate recorded in
Section~\ref{ssec:total-error} of the supplementary materials, and its size is assessed by
independent prices and refinement (Section~\ref{ssec:calibration}).
\subsection{Spectral lower bound}
\label{ssec:spectral}
\begin{lemma}[L\'evy--Khintchine form of the real symbol]\label{lem:LK}
For $Y\in(0,2)\setminus\{1\}$, $a>0$ and $\xi\in\R$,
\[
 C\Gamma(-Y)\bigl[a^Y-\Real(a+i\xi)^Y\bigr]
 =C\,I_a(\xi),\qquad
 I_a(\xi):=\int_0^\infty(1-\cos\xi x)\,e^{-ax}x^{-1-Y}dx\ \ge0 .
\]
Consequently $\Real A_\mu(\xi)=\rf+C\,[I_G(\xi)+I_M(\xi)]\ge \rf$ for all $\xi$,
without any restriction $Y>1$.
\end{lemma}
\begin{proof}
For $Y\in(1,2)$: $\int_0^\infty(e^{-zx}-1+zx)x^{-1-Y}dx=\Gamma(-Y)z^Y$,
$\Real z\ge0$; apply with $z=a+i\xi$ and $z=a$, subtract, use
$\Gamma(1-Y)=-Y\Gamma(-Y)$ and take real parts.  For $Y\in(0,1)$ use
$\int_0^\infty(e^{-zx}-1)x^{-1-Y}dx=\Gamma(-Y)z^Y$ in the same way; the
sign of $\Gamma(-Y)$ cancels against that of the bracket.
\end{proof}

\begin{theorem}[Rigorous spectral lower bound]\label{thm:spectral}
For exact assembly and every $N$,
\begin{equation}
 \lambda_{\min}(S_I)\ \ge\ \smin(\mu):=\min_{\vartheta}\sym{s}_\mu(\vartheta)
 \ \ge\ \rf\,m_{\min},
 \qquad S_I\succeq \rf\,M_I .
 \label{eq:spectral}
\end{equation}
For the truncated assembly with $K$ aliases the same bounds hold for
$\Tn{N_h}{\sym{s}_{\mu,K}}$ with $\min\sym{s}_{\mu,K}$ in place of $\smin$, and both
matrices satisfy $\lambda_{\min}\ge\alpha_{\rm lb}(\mu):=\min\sym{s}_{\mu,K}-\varepsilon_K(\mu)$.
All bounds are parameter-explicit and valid for $Y\in(0,2)\setminus\{1\}$.
\end{theorem}
\begin{proof}
By Lemma~\ref{lem:LK},
\[
\sym{s}_\mu(\vartheta)=\sum_k|\widehat N_p(\vartheta_k)|^2\Real A_\mu(\vartheta_k/h)\ge \rf\,\sym{m}(\vartheta),
\]
and Lemma~\ref{lem:toeplitz} applies.  For the truncated symbol use
$|\sym{s}_\mu-\sym{s}_{\mu,K}|\le\varepsilon_K$ (Lemma~\ref{lem:tail}) and
$\norm{\Tn{N_h}{\sym{s}_\mu}-\Tn{N_h}{\sym{s}_{\mu,K}}}_2\le\varepsilon_K$.
\end{proof}

\begin{remark}[Assembled versus exact sections]\label{rem:compression}
The matrices actually assembled are the sections built from the $L$-point DFT coefficients of
$\sym{s}_{\mu,K}$ (Section~\ref{ssec:assembly}).  Up to the diagonal unitary similarity
$\operatorname{diag}((-1)^k)$, such a section is the leading $N_h\times N_h$ compression of the
$L\times L$ circulant whose eigenvalues are the samples $\sym{s}_{\mu,K}(\vartheta_j)$, so for the
assembled matrix $\lambda_{\min}(S_I)\ge\min_j\sym{s}_{\mu,K}(\vartheta_j)$ holds \emph{exactly},
and likewise $\norm{\Tn{N_h}{f}-\Tn{N_h}{g}}_2\le\max_j|f(\vartheta_j)-g(\vartheta_j)|$ for two
symbols sampled on the same grid.  The derivative margin of the grid minimum and the tail
$\varepsilon_K$ are needed only to relate the assembled system to the exact Galerkin matrix; the
implementation keeps both, so that $\alpha_{\rm lb}$ is valid for either matrix.
\end{remark}
Over 32 box vertices and 256 Sobol points at $N=64$, $256$ and $1024$,
$\lambda_{\min}(S_I)/\smin(\mu)\in[1.003,17.8]$ and $\lambda_{\min}(S_I)/(\rf m_{\min})\in[7.5,133]$
(Figure~\ref{fig:spectral-bound}), the ratio closest to one being attained at the worst vertex of
the earlier positivity audit, $\mu=(0.05,5,9,1.1,0.01)$, where $\lambda_{\min}=0.01038$ and
$\smin=\rf=0.0100$ hold.
$\smin(\mu)$ is the minimum of an explicit smooth even function of one variable, evaluated on a
graded grid ($2^{14}$ uniform points plus a geometric refinement on the scale $\min(G,M)h$ near
$\vartheta=0$, where the fractional part varies fastest) and shifted outward by a cell-wise
derivative margin $\tfrac12 w\,\sup_{\rm cell}|\sym{s}'_{\mu}|$ with the derivative from
$\partial_\vartheta(\delta+i\vartheta)^Y=iY(\delta+i\vartheta)^{Y-1}$; the only step that is not an
enclosure is that $\sup_{\rm cell}|\sym{s}'_\mu|$ is taken as the larger endpoint value (a validated
implementation would use interval arithmetic on each cell).  On the box the minimum is attained at
$\vartheta=0$, where $\sym{s}_\mu(0)=\rf$ exactly, so in practice $\alpha_{\rm lb}=\rf-\varepsilon_K(\mu)$
and the margin is of the order $10^{-8}$.
Since the certificate needs only $S_I\succ0$, all algebraic results extend to $Y\in(0,1)$
(Section~\ref{ssec:e0-initial}); only the continuum coercivity in $H^{Y/2}_{M,G}$ keeps the
hypothesis $Y>1$.
\section{Residual certification and tolerance-controlled correction}
\label{sec:certificate}\label{sec:pcc}
\subsection{Energy norm and Crank--Nicolson bound}
\label{ssec:riesz}\label{ssec:cn-certificate}
For a fixed query, assume $\Mh\succ0$ and
$\Sh=(\Ah+\Ah^\top)/2\succ0$ on the free space. Define
\begin{equation}
 \|v\|_{\Sh}^2=v^\top\Sh v,\qquad
 \|r\|_{\Sh^{-1}}^2=r^\top\Sh^{-1}r
 =\sup_{v\ne0}\frac{(r^\top v)^2}{v^\top\Sh v}.
 \label{eq:Sh}
\end{equation}
Positivity of $S_I$ holds by Theorem~\ref{thm:spectral} with the explicit
constant $\alpha_{\rm lb}(\mu)=\min\sym{s}_{\mu,K}-\varepsilon_K(\mu)>0$.

For an arbitrary candidate trajectory $\{\cc_\theta^n\}_{n=0}^{N_T}$ define the
Crank--Nicolson Galerkin residual
\begin{equation}
 \boxed{
 \rr_\theta^{n+1/2}
 :=\Mh\frac{\cc_\theta^{n+1}-\cc_\theta^n}{\dt}
  +\Ah(\mu)\overline\cc_\theta^{\,n+1/2}
  -\overline\ff^{\,n+1/2},}
 \label{eq:cn-residual}
\end{equation}
and the coefficient error $\ee^n=\cc_h^n-\cc_\theta^n$.  Subtracting
\eqref{eq:cn-residual} from \eqref{eq:cn-reference},
\begin{equation}
 \Mh\frac{\ee^{n+1}-\ee^n}{\dt}+\Ah(\mu)\overline\ee^{\,n+1/2}
 =-\rr_\theta^{n+1/2}.
 \label{eq:discrete-error-equation}
\end{equation}

\begin{theorem}[Crank--Nicolson residual certificate]
\label{thm:fully-discrete-certificate}
Let $\Mh$ and $\Sh$ be symmetric positive definite, and let reference and
candidate use identical matrices, forcing and prescribed coefficients.
For every $m=1,\ldots,N_T$,
\begin{equation}
 \|\ee^m\|_{\Mh}^2+
 \dt\sum_{n=0}^{m-1}\|\overline\ee^{n+1/2}\|_{\Sh}^2
 \le\|\ee^0\|_{\Mh}^2+
 \dt\sum_{n=0}^{m-1}\|\rr_\theta^{n+1/2}\|_{\Sh^{-1}}^2.
 \label{eq:fully-discrete-certificate}
\end{equation}
Consequently, with
\begin{equation}
 \etaSharp^2=\dt\sum_{n=0}^{N_T-1}
 (\rr_\theta^{n+1/2})^\top\Sh^{-1}\rr_\theta^{n+1/2},
 \qquad \eta_{\rm aug,sh}=(\beta_0^2+\etaSharp^2)^{1/2},
 \label{eq:sharp-estimator}
\end{equation}
where $\beta_0\ge\|\ee^0\|_{\Mh}$, one has
\begin{equation}
 \max_{0\le m\le N_T}\|V_h^m-V_{\theta,h}^m\|_{L^2(\Omega)}
 \le\eta_{\rm aug,sh}.
 \label{eq:sharp-L2}
\end{equation}
\end{theorem}
\begin{proof}
Multiply \eqref{eq:discrete-error-equation} on the left by $(\overline\ee^{\,n+1/2})^\top$: the
mass term equals $(\norm{\ee^{n+1}}_{\Mh}^2-\norm{\ee^n}_{\Mh}^2)/(2\dt)$ because $\Mh$ is
symmetric, the operator term equals $\norm{\overline\ee^{\,n+1/2}}_{\Sh}^2$ by \eqref{eq:Sh}, and
the residual term is bounded by Young's inequality,
\[
|(\overline\ee^{\,n+1/2})^\top\rr_\theta^{n+1/2}|\le\tfrac12\norm{\overline\ee^{\,n+1/2}}_{\Sh}^2
+\tfrac12\norm{\rr_\theta^{n+1/2}}_{\Sh^{-1}}^2 .
\]
Multiply by $2\dt$ and sum over
$n=0,\dots,m-1$; the mass terms telescope and give \eqref{eq:fully-discrete-certificate}.  The
right-hand side increases with $m$; its final value and the shared-lifting cancellation of
Section~\ref{ssec:reference-time} prove \eqref{eq:sharp-L2}.
\end{proof}
With matching initial vectors the estimator vanishes precisely when the trajectories coincide,
since the symmetric part of $B_+=\Mh+\dt\Ah/2$ is positive definite and the CN recurrence is
uniquely solvable.  If the residual is evaluated with different matrices or forcing, the required
residual is $\rho_{\rm ref}=\rho_{\rm ev}+(A_h-A_{\rm ev})\overline c_\theta+(f_{\rm ev}-f_{\rm ref})$,
and the discrepancies require allowances in the same dual norm (a computed reference also carries a
solver-error allowance).

\subsection{Upper bounds from incomplete Riesz solves}
\label{ssec:guaranteed}
In practice the Riesz problem $\Sh v=\rr$ is solved iteratively, so the online method needs an
upper bound that remains valid before the iteration has converged.  It uses the lower bound
\begin{equation}
  0<\alpha_{\rm lb}(\mu)\le \lambda_{\min}(\Sh(\mu)),
  \qquad \alpha_{\rm lb}(\mu):=\min_\vartheta\sym{s}_{\mu,K}(\vartheta)-\varepsilon_K(\mu),
  \label{eq:alpha-lb}
\end{equation}
which Theorem~\ref{thm:spectral} supplies in closed form; no eigenvalue computation enters.
\begin{lemma}[Computable two-sided bracket for the dual norm]
\label{lem:guaranteed}
Let $\Sh\succ0$, $\rr\in\R^{N_h}$ and $0<\alpha_{\rm lb}\le\lambda_{\min}(\Sh)$.  For an arbitrary
$\vv\in\R^{N_h}$ set $\sres:=\rr-\Sh\vv$ and
\begin{equation}
  \mathfrak l(\vv):=2\vv^\top\rr-\vv^\top\Sh\vv,
  \qquad
  \mathfrak u(\vv):=\mathfrak l(\vv)+\frac{\norm{\sres}_2^2}{\alpha_{\rm lb}}.
  \label{eq:riesz-bracket-def}
\end{equation}
Then $\mathfrak l(\vv)\le\rr^\top\Sh^{-1}\rr\le\mathfrak u(\vv)$ for every $\vv$, with width
$\norm{\sres}_2^2/\alpha_{\rm lb}$, which closes as the Riesz residual tends to zero.
\end{lemma}
\begin{proof}
With $\zz=\Sh^{-1}\rr$, $\rr^\top\Sh^{-1}\rr-\mathfrak l(\vv)=(\vv-\zz)^\top\Sh(\vv-\zz)
=\sres^\top\Sh^{-1}\sres\ge0$, and $\sres^\top\Sh^{-1}\sres\le\alpha_{\rm lb}^{-1}\norm{\sres}_2^2$
since $\Sh\succeq\alpha_{\rm lb}I$.
\end{proof}
\begin{corollary}[Augmented upper certificate]
\label{cor:guaranteed}
Under the hypotheses of Theorem~\ref{thm:fully-discrete-certificate} and
$0<\alpha_{\rm lb}\le\lambda_{\min}(\Sh)$, let $v^{n+1/2}$ be any approximate Riesz vectors and
\begin{equation}
 \etaGuar^2=\dt\sum_{n=0}^{N_T-1}\mathfrak u(v^{n+1/2}),
 \qquad \eta_{\rm aug}=(\beta_0^2+\etaGuar^2)^{1/2}.
 \label{eq:guaranteed-estimator}
\end{equation}
Then $\max_m\|\ee^m\|_{\Mh}\le\eta_{\rm aug}$, independently of the iterative stopping error,
and $\etaGuar\to\etaSharp$ as the Riesz residuals vanish.
\end{corollary}
The zero iterate gives the upper bound $\|\rr\|_2^2/\alpha_{\rm lb}$; convergence gives the exact
Riesz value.
\begin{remark}[Arithmetic scope]
\label{rem:arithmetic-scope}
The theorems and brackets are exact-arithmetic statements for the assembled matrices.
The only quantities not enclosed are the floating-point residual evaluation and the
accumulation: if a computed residual $\widehat\rho_n$ differs from the required residual by at
most $d_n$ in the $S_h^{-1}$ norm and $u_n$ is an upper bound for
$\widehat\rho_n^\top S_h^{-1}\widehat\rho_n$, then
$\widehat\eta^2=\beta_0^2+\sum_n\dt\,(\sqrt{u_n}+d_n)^2$ is valid.  The same allowance
structure absorbs the affine-operator error of Theorem~\ref{thm:eim-cert}.
\end{remark}
\subsection{A banded comparison norm}
\label{ssec:comparison}
The dual norm $\rr^\top S_I^{-1}\rr$ requires a solve with the dense Toeplitz matrix $S_I$ at
every time step.  Lemma~\ref{lem:toeplitz} suggests replacing $S_I$ by a banded Toeplitz matrix
whose symbol lies below $\sym{s}_\mu$: the dual norm in the banded matrix is then an upper bound,
with an explicit loss factor.  The natural minorant is built from the two symbols that are
available in closed form, the mass symbol and the stiffness symbol $\sym{g}_1$, since by
Lemma~\ref{lem:LK} the fractional part of $\sym{s}_\mu$ is nonnegative and vanishes
quadratically at $\vartheta=0$ exactly as $\sym{g}_1$ does.
\begin{theorem}[Comparison norm]\label{thm:comparison}
Fix $\mu$ and define
\[
 \beta_\mu:=\inf_{\vartheta\ne0}\frac{\sym{s}_\mu(\vartheta)-\rf\,\sym{m}(\vartheta)}{\sym{g}_1(\vartheta)}>0,
 \qquad
 \Thmu:=\sup_\vartheta\frac{\sym{s}_\mu(\vartheta)}{\sym{x}_\mu(\vartheta)},
 \qquad
 \sym{x}_\mu:=\rf\,\sym{m}+\beta_\mu \sym{g}_1,
\]
and the banded matrix
\[
 X_I:=\Tn{N_h}{\sym{x}_\mu}=\rf\,M_I+\beta_\mu h^2R^{(1)}_I.
\]
$X_I$ is symmetric positive definite and banded with bandwidth $p-1$, and
for every $\rho\in\R^{N_h}$
\begin{equation}
 \rho^\top S_I^{-1}\rho\ \le\ \rho^\top X_I^{-1}\rho\ \le\ \Thmu\,\rho^\top S_I^{-1}\rho .
 \label{eq:comparison}
\end{equation}
\end{theorem}
\begin{proof}
$\sym{s}_\mu-\rf\,\sym{m}=C\sum_k|\widehat N_p(\vartheta_k)|^2[I_G+I_M](\vartheta_k/h)\ge0$ vanishes
only at $\vartheta=0$, where both it and $\sym{g}_1$ are $\asymp\vartheta^2$
($I_a(\xi)\sim\tfrac12 I_a''(0)\xi^2$ with $I_a''(0)=\Gamma(2-Y)a^{Y-2}>0$), so
$\beta_\mu>0$.  Lemma~\ref{lem:toeplitz} gives $S_I\succeq X_I\succeq\Thmu^{-1}S_I$;
invert in the Loewner order.
\end{proof}

\begin{corollary}[Cheap certificate]\label{cor:cheap}
In Theorem~\ref{thm:fully-discrete-certificate} and Proposition~\ref{prop:two-stage} replace
$\norm{\rr^{n+1/2}}_{\Sh^{-1}}^2$ by $\rho^\top X_I^{-1}\rho$.  The resulting
$\eta^{X}_{\rm aug}=(\beta_0^2+(\eta^X)^2)^{1/2}$ with
$(\eta^X)^2=\dt\sum_n(\rr^{n+1/2})^\top X_I^{-1}\rr^{n+1/2}$
is a rigorous upper bound with effectivity at most
$\sqrt{\Thmu}$ times that of the exact estimator $\eta^S_{\rm aug}$, and costs one banded
Cholesky factorisation plus $N_T$ banded solves ($O(N_TN_h)$ operations).
No Riesz solve, no spectral estimate and no bracket are needed.
\end{corollary}

From here on $\eta^S$ denotes the exact estimator $\etaGuar$ at the converged Riesz iterate
($\eta^S=\etaSharp$ up to the bracket width), $\eta^X$ the banded estimator of
Corollary~\ref{cor:cheap}, and $\eta^S_{\rm aug},\eta^X_{\rm aug}$ their augmented versions
$(\beta_0^2+\eta^2)^{1/2}$.
\begin{proposition}[Two-stage rule]\label{prop:twostage}
Let $\eta^X$ be the cheap residual estimator and $\eta^S$ the exact one, with the common
initial term $\beta_0$.  Then
$\eta^S\le\eta^X\le\sqrt{\Thmu}\,\eta^S$.  Hence, for a tolerance $\varepsilon$:
if $(\beta_0^2+(\eta^X)^2)^{1/2}\le\varepsilon$ accept (\textsc{banded accept});
if $(\beta_0^2+(\eta^X)^2/\Thmu)^{1/2}>\varepsilon$ reject without further work
(\textsc{cheap reject}: the exact bound cannot pass);
otherwise compute $\eta^S$ by the bracket of Lemma~\ref{lem:guaranteed} with
$\alpha_{\rm lb}=\min\sym{s}_{\mu,K}-\varepsilon_K(\mu)$ and decide
(\textsc{exact accept} / \textsc{exact reject}).
\end{proposition}

The two-stage rule never evaluates the exact bound when the cheap one already decides, and never
accepts on the basis of a lower bound.  For the third case the Riesz systems $S_Iv=\rho$ are
solved by conjugate gradients with the circulant preconditioner whose eigenvalues are the symbol
samples $\sym{s}_\mu(2\pi j/N_h)$, $j=0,\dots,N_h-1$, read off the same grid on which the operator
was assembled.  Since $\sym{s}_\mu$ is continuous and positive, the preconditioned spectrum
clusters at $1$ \cite{ChanStrang1989}; Table~\ref{tab:pcg-v9} reports $8$--$12$ iterations to a
relative residual of $10^{-12}$ for $N=64,\dots,4096$, against $42$--$1915$ for plain conjugate
gradients, and the bracket width at the final iterate is below $10^{-18}$ relative.  Both
quantities that enter the certificate per time step, the residual $\rr^{n+1/2}$ (two FFTs of
length $2N_h$ for $A_I$, a banded product for $M_I$) and the banded dual norm (one banded
Cholesky factorisation of $X_I$ per query, one banded solve per step), cost $O(N_h\log N_h)$ and
$O(N_h)$ respectively; the certificate is therefore cheaper than one Krylov step of the direct
solve, which is the one-dimensional ceiling discussed in Section~\ref{ssec:cost}.

When the candidate comes with the $N$-independent reduced estimator of Section~\ref{ssec:eim}, the
rule gains a first stage that never touches the full space:
\begin{proposition}[Three-stage rule]\label{prop:threestage}
Let $\bar X=\Tn{N_h}{\bar{\sym{x}}}$ be the comparison matrix of a cell containing $\mu$
(Remark~\ref{rem:reduced-online}), so that $\bar X\preceq X_I\preceq S_I$, and let
$\bar\Theta_\mu:=\sup_\vartheta\sym{s}_\mu(\vartheta)/\bar{\sym{x}}(\vartheta)$, so that
$S_I\preceq\bar\Theta_\mu\bar X$ by Lemma~\ref{lem:toeplitz}.  Let $\rho_Q^{n+1/2}$ be the residual of a
candidate evaluated with the affine operator $A_Q$, $u_n\ge(\rho_Q^{n+1/2})^\top\bar X^{-1}\rho_Q^{n+1/2}$
the per-step values of the reduced online estimator (including its rounding term), and $d_n$ the
allowance of Theorem~\ref{thm:eim-cert} with $\min_\vartheta\bar{\sym{x}}$ in place of $\alpha_{\rm lb}$.  Then
\begin{align*}
 \eta^{\bar X}_{\rm aug}&:=\Bigl(\beta_0^2+\dt\sum_n(\sqrt{u_n}+d_n)^2\Bigr)^{1/2}\ \ge\ \eta^S_{\rm aug}\\
 &\ \ge\ \Bigl(\beta_0^2+\dt\sum_n\max\bigl\{0,\ \sqrt{\ell_n/\bar\Theta_\mu}-d_n\bigr\}^2\Bigr)^{1/2}=:\lambda^{\bar X}_{\rm aug},
\end{align*}
where $\ell_n\le(\rho_Q^{n+1/2})^\top\bar X^{-1}\rho_Q^{n+1/2}$ is the per-step value minus its rounding
term.  Hence, for a tolerance $\varepsilon$: if $\eta^{\bar X}_{\rm aug}\le\varepsilon$ accept
(\textsc{reduced accept}); if $\lambda^{\bar X}_{\rm aug}>\varepsilon$ reject (\textsc{cheap reject}:
no full-order certificate can pass); otherwise apply Proposition~\ref{prop:twostage} at full order
with the exact operator $\Ah$ (\textsc{banded accept}, \textsc{exact accept}, \textsc{exact reject}).
Only the first two stages are independent of $N$.
\end{proposition}
\begin{proof}
$\rho^{n+1/2}=\rho_Q^{n+1/2}+E^\top(\bar A_h-\bar A_Q)\overline{\bar\cc}_\theta$ with
$\norm{E^\top(\bar A_h-\bar A_Q)\overline{\bar\cc}_\theta}_{S_I^{-1}}\le d_n$ (Theorem~\ref{thm:eim-cert},
$S_I\succeq\bar X$); the upper bound follows from $S_I^{-1}\preceq\bar X^{-1}$ and the lower one from
$S_I^{-1}\succeq\bar\Theta_\mu^{-1}\bar X^{-1}$ (Loewner order of Lemma~\ref{lem:toeplitz}); then
Theorem~\ref{thm:fully-discrete-certificate}.
\end{proof}

The lower bound $\lambda^{\bar X}_{\rm aug}$ is tight in practice (within $7\%$ of $\eta^S_{\rm aug}$ on the
reduced candidates of Section~\ref{ssec:e1}), so the cheap rejection settles every clearly
inadequate candidate, and the full-order stages are reached only for candidates whose bound lies
within a factor $\sqrt{\bar\Theta_\mu}\approx4$ of the tolerance.
\subsection{Efficiency and spot-price bounds}
\label{ssec:efficiency}\label{ssec:greeks-cert}
For $\Sh\succ0$, define
\begin{equation}
 \kappa_h=\|\Sh^{-1/2}\Ah\Sh^{-1/2}\|_2\ge1,
 \qquad
 |||\ee|||^2=\dt\sum_n
 \left(\|\Mh\partial\ee^{n+1/2}\|_{\Sh^{-1}}^2
       +\|\overline\ee^{n+1/2}\|_{\Sh}^2\right),
 \label{eq:kappa-def}
\end{equation}
where $\partial\ee^{n+1/2}=(\ee^{n+1}-\ee^n)/\dt$.
\begin{theorem}[Matrix-dependent efficiency]
\label{thm:efficiency}
For an assembled query with $\Sh\succ0$,
\begin{equation}
 \etaSharp^2\le2\max\{1,\kappa_h^2\}\,|||\ee|||^2.
 \label{eq:efficiency}
\end{equation}
Under exact assembly, $\kappa_h\le C_{\rm cont}/c_1$, uniformly in the
mesh and parameter family.
\end{theorem}
\begin{proof}
The error equation gives
$\|\rr\|_{\Sh^{-1}}\le\|\Mh\partial\ee\|_{\Sh^{-1}}
+\kappa_h\|\overline\ee\|_{\Sh}$. Square, use
$(a+b)^2\le2a^2+2b^2$, and sum. The exact-assembly estimate follows from
continuity and coercivity; see Section~\ref{ssec:exact-assembly} of the supplementary materials.
\end{proof}
The efficiency norm contains a time-derivative term and differs from the maximum-in-time mass
error of the reliability estimate (for zero initial error, $|||\ee|||^2\le(3+2\kappa_h^2)\etaSharp^2$);
effectivity is measured as the augmented bound divided by $\max_m\|\ee^m\|_{\Mh}$.
For $\Phi_I(x)$ the free-basis evaluation vector, Cauchy--Schwarz in the mass pairing gives the
discrete spot-price bound
\begin{equation}
 |V_h^m(x)-V_{\theta,h}^m(x)|\le\beta_I(x)\eta_{\rm aug},\qquad
 \beta_I(x)=\sqrt{\Phi_I(x)^\top M_I^{-1}\Phi_I(x)},
 \label{eq:gamma-infty}
\end{equation}
which converts a trajectory tolerance into a price tolerance at specified spots; derivative
bounds and the continuum error separation are recorded in Section~\ref{ssec:derivatives} of the
supplementary materials.

\subsection{Tolerance-controlled prediction and correction}
\label{ssec:algorithm}\label{ssec:two-stage}\label{ssec:correction}
\begin{figure}[!htbp]
\centering
\begin{minipage}{0.96\linewidth}\small
\textbf{Tolerance-checked predict--certify--correct rule.}
\begin{enumerate}[leftmargin=2em,itemsep=1pt]
\item Construct the initial target, predict the coefficient trajectory, overwrite the prescribed
coefficients by the specified lifting at every time level, and compute
$\beta_0=\|\zz_\theta^0-\zz_h^0\|_{M_I}$.  Check that the time grid matches the uniform CN scheme
and form the free residuals with the same operator, stored forcing and lifting as the reference.
\item If the candidate comes from the reduced model with the affine operator, evaluate the
$N$-independent estimator $\eta^{\bar X}_{\rm aug}$ and its lower companion $\lambda^{\bar X}_{\rm aug}$
(Proposition~\ref{prop:threestage}): return \textsc{reduced accept} if the former is at most
$\varepsilon$, go to step 4 if the latter exceeds $\varepsilon$, and otherwise proceed at full order.
\item Evaluate $\alpha_{\rm lb}(\mu)$, $\beta_\mu$, $X_I$ and $\Thmu$ from the symbol
(Theorems~\ref{thm:spectral} and \ref{thm:comparison}) and the cheap augmented bound
$\eta^X_{\rm aug}=(\beta_0^2+(\eta^X)^2)^{1/2}$ (banded Cholesky).  If it is at most $\varepsilon$,
return \textsc{banded accept}.  If $(\beta_0^2+(\eta^X)^2/\Thmu)^{1/2}>\varepsilon$, no exact bound
can pass (\textsc{cheap reject}): go to step 4.  Otherwise compute $\eta^S_{\rm aug}$ by the
preconditioned bracket; if it passes, return \textsc{exact accept}, else \textsc{exact reject}.
\item Correct by a free CN march from $\zz_h^0$, warm-started by the candidate, and recertify the
corrected trajectory: \textsc{corrected accept} if its bound passes; otherwise a direct fallback,
recertified (\textsc{qualified fallback}), or \textsc{over budget} / \textsc{numerical failure}.  A
certified output is returned only when the recomputed bound passes.
\end{enumerate}
\end{minipage}
\caption{Tolerance-checked online rule (the three-stage rule of Proposition~\ref{prop:threestage}
followed by correction). The status names are those used in the tables of Section~\ref{sec:numerics}.}
\label{alg:pcc}
\end{figure}
The bracket is evaluated at the converged PCG iterate (relative residual $10^{-12}$, at most
$12$ iterations in all reported runs, Table~\ref{tab:pcg-v9}); correction uses circulant-preconditioned
GMRES with relative tolerance $10^{-13}$, and tolerances are reported as a sweep
(Section~\ref{ssec:calibration}).
\begin{proposition}[Guaranteed acceptance]
\label{prop:two-stage}
Assume the hypotheses of Theorem~\ref{thm:fully-discrete-certificate},
exact residual evaluation, and $\beta_0\ge\|\ee^0\|_{\Mh}$. Let $\eta_{\rm aug}$ be
$\eta^{\bar X}_{\rm aug}$, $\eta^X_{\rm aug}$ or $(\beta_0^2+\etaGuar^2)^{1/2}$ with $\alpha_{\rm lb}$ from \eqref{eq:alpha-lb}.
If the rule in Figure~\ref{alg:pcc} returns
\textsc{reduced accept}, \textsc{banded accept}, \textsc{exact accept} or \textsc{corrected accept} with $\eta_{\rm aug}\le\varepsilon$, then
\begin{equation}
  \max_{1\le m\le N_T}\norm{V_h^m-V_{\theta,h}^m}_{L^2(\Omega)}
  \le \eta_{\rm aug}\le\varepsilon .
  \label{eq:soundness}
\end{equation}
The conclusion is independent of the PCG stopping error, by Lemma~\ref{lem:guaranteed} for
$\etaGuar$, by Corollary~\ref{cor:cheap} for $\eta^X$ and by Proposition~\ref{prop:threestage} for $\eta^{\bar X}$;
it applies to a corrected trajectory or a direct fallback whose recomputed augmented bound passes.
Bounded evaluation discrepancies additionally require the allowances of Remark~\ref{rem:arithmetic-scope}.
\end{proposition}
\begin{proof}
At every time step the computable upper bracket bounds $\rr^\top\Sh^{-1}\rr$; sum and apply
Theorem~\ref{thm:fully-discrete-certificate} with the initial term bounded by $\beta_0^2$.
\end{proof}
The correction $\dd^n$ with $\cc^n_{\rm corr}=\cc_\theta^n+\dd^n$ that recovers the reference
solves the error equation \eqref{eq:discrete-error-equation} with $\dd^0=\ee^0$, i.e.\ with
$B_\pm=\Mh\pm\tfrac{\dt}{2}\Ah$,
\begin{equation}
  B_+\dd^{n+1}=B_-\dd^n-\dt\,\rr_\theta^{n+1/2},\qquad \dd^0_I=\zz_h^0-\zz_\theta^0;
  \label{eq:correction-step}
\end{equation}
solving it exactly is equivalent to recomputing the reference, so the practical aim is a small
number of Krylov iterations per step warm-started with the candidate coefficients, and an
approximate correction is accepted only after its augmented bound is recomputed.  A conditional
Richardson estimate is given in Section~\ref{ssec:correction-extension} of the supplementary
materials; the implementation uses GMRES.

\section{Affine structure and surrogate candidates}
\label{sec:affine}\label{sec:operator}
The certificate of Section~\ref{sec:certificate} applies to any trajectory in the spline space.
This section describes the three families of candidates evaluated in Section~\ref{sec:numerics},
a projection-based reduced model, a neural reduced model and a neural operator, and the
structure that makes the first two certifiable at a cost independent of $N$: the operator
family $\{\Ah(\mu)\}$ is, to a stated accuracy, affine in a small number of
parameter-independent Toeplitz matrices.
\subsection{Affine decomposition of the fractional blocks}
\label{ssec:eim}
The operator depends on $\mu$ through the scalars $\rf-\omega_S$, $\lambda_0$, $C\Gamma(-Y)h^{-Y}$
and, non-affinely, through $(G,M,Y)$ in the fractional symbol.  The reduced-basis literature
handles such dependence by the empirical interpolation method
\cite{BarraultMadayNguyenPatera2004,GreplMadayNguyenPatera2007}; what is specific here is that
the interpolation acts on the symbol, a function of one variable, rather than on the matrices,
and that the interpolation error transfers to the operator norm through
Lemma~\ref{lem:toeplitz} without any constant.
Write
\begin{equation}
\begin{gathered}
 \sym{a}_\mu=-(\rf-\omega_S)\tfrac{i}{h}\sym{d}_1+\lambda_0\,\sym{m}-C\Gamma(-Y)h^{-Y}\Phi_{\nu},\\
 \Phi_\nu(\vartheta):=\sum_k|\widehat N_p(\vartheta_k)|^2\bigl[(Gh+i\vartheta_k)^Y+(Mh-i\vartheta_k)^Y\bigr],
\end{gathered}
 \label{eq:affine-split}
\end{equation}
with $\nu=(G,M,Y)$.  The first two terms are already affine in $\mu$ with
fixed matrices $D_h$, $\Mh$.  What has to be interpolated is the
three-parameter family of periodic functions $\{\Phi_\nu\}$, not the
matrix family: because $\Tn{n}{\cdot}$ is linear, an interpolant
$\Phi_\nu\approx\sum_{q=1}^Q\ell_q(\nu)\Phi_{\nu_q}$ gives
$A_Q(\mu):=-(\rf-\omega_S)D_h+\lambda_0\Mh-C\Gamma(-Y)h^{-Y}\sum_q\ell_q(\nu)\Tn{}{\Phi_{\nu_q}}$
with $Q+2$ terms and the same coefficient functions.
\begin{lemma}[Error transfer]\label{lem:eim}
Let $e_\nu:=\Phi_\nu-\sum_q\ell_q(\nu)\Phi_{\nu_q}$.  Then
\[
 \norm{\Ah(\mu)-A_Q(\mu)}_2\ \le\ \varepsilon_Q(\mu):=C|\Gamma(-Y)|h^{-Y}\norm{e_\nu}_\infty,
\]
and $\norm{e_\nu}_\infty\le\max_{j}|e_\nu(\vartheta_j)|+\tfrac{\pi}{L}\norm{e_\nu'}_\infty$
on an $L$-point grid, with $e_\nu'$ obtained from the explicit derivative
$\partial_\vartheta(\delta+i\vartheta)^Y=iY(\delta+i\vartheta)^{Y-1}$ and the
same tail bound.  For the assembled matrices (DFT sections on the same grid) the grid maximum
alone is a rigorous bound, by Remark~\ref{rem:compression}.  Evaluating $\varepsilon_Q(\mu)$ at
a new query costs $O(LK)$ scalar operations and no linear algebra.
\end{lemma}
\begin{proof}
Lemma~\ref{lem:toeplitz}, second inequality, and the mean value theorem.
\end{proof}

\begin{figure}[!htbp]
\centering
\includegraphics[width=0.8\linewidth]{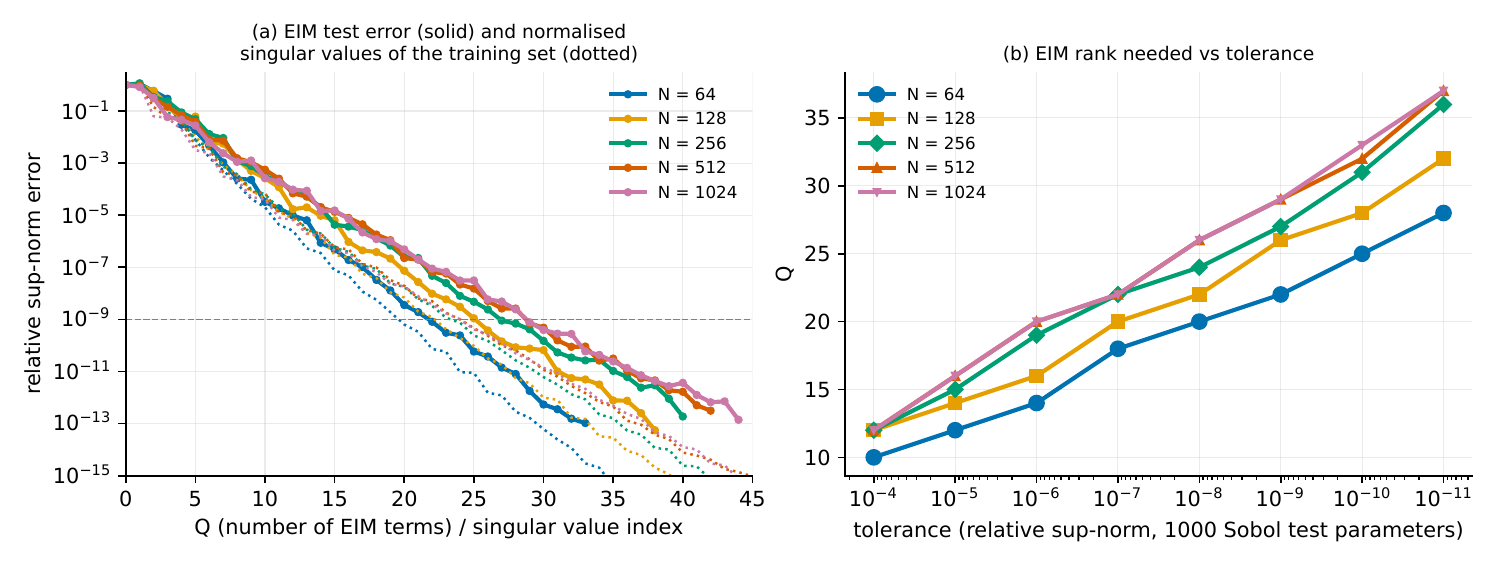}
\caption{Empirical interpolation of the fractional symbol family.  (a) Relative sup-norm error on
$1000$ Sobol test parameters (solid) and normalised singular values of the $16^3$ training snapshots
(dotted), $N=64$--$1024$.  (b) Rank $Q$ needed for a given test tolerance; the curves for $N=512$ and
$N=1024$ coincide.  The operational $Q$ additionally satisfies the operator-norm check
(one extra term at $N=64$).}
\label{fig:eim}
\end{figure}

Greedy empirical interpolation of $\{\Phi_\nu\}$ on a $16\times16\times16$ training grid in
$(\log G,\log M,Y)$ selects magic points $\vartheta_1,\dots,\vartheta_Q$ and parameters
$\nu_1,\dots,\nu_Q$ such that $\ell(\nu)$ is the solution of a $Q\times Q$ triangular system; the
family is analytic in $(\delta,Y)$ with its nearest singularity at $\delta=0$, so the
interpolation error decays geometrically in $Q$ at a rate set by $G_{\max}/G_{\min}$ and
$M_{\max}/M_{\min}$ rather than by $h$.  Figure~\ref{fig:eim} shows the test error on $1000$ Sobol parameters and the singular values
of the training snapshots: $Q=23$, $26$, $27$, $29$, $29$ terms reach a relative sup-norm error of
$10^{-9}$ at $N=64$, $128$, $256$, $512$, $1024$, and the rank saturates in $h$.  For the operator the relevant number is $\varepsilon_Q$ relative to
$\norm{\Ah}_2$: at most $1.7\times10^{-9}$ ($N=64$) and $5.7\times10^{-10}$ ($N=512$) with these $Q$,
with the observed operator-norm error never exceeding $0.93\,\varepsilon_Q$.  The operational rank
is chosen per target tolerance by the condition $\varepsilon_Q\norm{\bar\cc}_2\sqrt{T/\alpha_{\rm lb}}
\le\varepsilon/10$ of Theorem~\ref{thm:eim-cert}: $Q=36$ at $\varepsilon_N$ and $Q=37$ at
$\varepsilon_N/10$ for the production configuration of Section~\ref{ssec:calibration}.
The error $\varepsilon_Q(\mu)$ is evaluated, not estimated: it is the sup-norm of the
interpolation error of the exact symbol at the query with the grid-to-continuum margin of
Lemma~\ref{lem:eim}.  Evaluated on the assembly grid it costs $O(LK)$ operations, which grows with
$N$; the implementation instead evaluates it on a fixed set of $161$ nodes (the graded grid of
Section~\ref{ssec:spectral}, coarsened) with the cell-wise derivative margin, which bounds the
continuous sup-norm and hence, by Remark~\ref{rem:compression}, the assembled operator difference,
at a cost of about $1.6$~ms independent of $N$ and at most $2.1$ times the value on the full grid.
Outside the training box of the interpolant the same evaluation remains rigorous and simply grows
(Section~\ref{ssec:ood}).
\begin{theorem}[Certificate with an affine surrogate operator]\label{thm:eim-cert}
Let the candidate trajectory be produced by any method (e.g.\ the reduced
system $V^\top A_Q(\mu)V$), and let $\rho_Q^{n+1/2}$ be its residual
evaluated with $A_Q$.  Then its residual with respect to the true assembled
system satisfies
\[
 \norm{\rr^{n+1/2}}_{S_I^{-1}}\ \le\ \norm{\rho_Q^{n+1/2}}_{S_I^{-1}}
 +\varepsilon_Q(\mu)\,\alpha_{\rm lb}^{-1/2}\,\norm{\overline{\bar\cc}_\theta^{\,n+1/2}}_2,
 \qquad \alpha_{\rm lb}=\smin(\mu),
\]
where $\bar\cc_\theta^n=Ez_\theta^n+b^n$ is the raw coefficient vector including the prescribed
entries (the surrogate operator also enters the lifting term of $f_I$), and
Theorem~\ref{thm:fully-discrete-certificate} holds with
$\eta^2=\beta_0^2+\dt\sum_n(\norm{\rho_Q}_{S^{-1}}+\varepsilon_Q\alpha_{\rm lb}^{-1/2}\norm{\overline{\bar\cc}_\theta}_2)^2$.
The same holds with $X_I$ in place of $S_I$ and $\min_\vartheta \sym{x}_\mu(\vartheta)\ (\ge \rf\,m_{\min})$ in place of $\alpha_{\rm lb}$.
\end{theorem}
\begin{proof}
$\rho=\rho_Q+E^\top(\bar A_h-\bar A_Q)\overline{\bar\cc}_\theta$ with the raw sections $\bar A_h,\bar A_Q$,
for which Lemma~\ref{lem:eim} gives the same $\varepsilon_Q$; bound the second term by
$\norm{(\bar A_h-\bar A_Q)\overline{\bar\cc}}_2\lambda_{\min}(S_I)^{-1/2}$ and use
Theorem~\ref{thm:spectral}.  This is exactly the allowance structure of
Remark~\ref{rem:arithmetic-scope}, with $d_n$ now explicit.
\end{proof}
The allowance is of size $\varepsilon_Q\norm{\bar\cc}_2/\sqrt{\alpha_{\rm lb}}$ per step and is
\emph{not} negligible for accurate candidates: with $\varepsilon_Q\approx10^{-9}\norm{\Ah}_2$,
$\norm{\Ah}_2\sim h^{-Y}$ and $\norm{\bar\cc}_2\sim \mathcal K\sqrt{N}\sqrt h$, its time-accumulated size at
$N=512$ with $Q=29$ is $6.8\times10^{-5}$ (median) to $4.3\times10^{-4}$ (maximum) in the $M_I$ norm on the
validation queries, comparable with the residual of a rank-$56$ reduced solution; its maximum drops to
$3.9\times10^{-6}$ with $Q=37$.  Two
consequences: the EIM tolerance must be chosen as a function of the target,
$\varepsilon_Q\norm{\bar\cc}_2\sqrt{T/\alpha_{\rm lb}}\le\varepsilon/10$, and the full-order
certificate should evaluate the residual with the exact operator $\Ah$ (two FFTs per step, no
allowance), reserving $A_Q$ for the reduced solve and for the $N$-independent online certificate.
\begin{remark}[Reduced online certificate]\label{rem:reduced-online}
With a fixed comparison norm $\bar X=\Tn{}{\bar{\sym{x}}}$, $\bar{\sym{x}}=\bar r\,\sym{m}+\bar\beta\,\sym{g}_1$,
chosen for a cell of parameter space such that $\bar r\le \rf$ and $\bar\beta\le\beta_\mu$ on the cell,
and the affine $A_Q$, the offline/online decomposition of the residual norm applies: the online cost of
$\eta^{\bar X}$ for a rank-$r$ candidate is $O(\varrho\,(Q+2)(2r+1))$ operations with $\varrho\le N_h$ retained rows
of the compressed $R$-factor, independent of $N_h$ once $\varrho$ saturates.  Since
$\bar{\sym{x}}\le\sym{x}_\mu\le\sym{s}_\mu$ on the cell, $\rho^\top S_I^{-1}\rho\le\rho^\top\bar X^{-1}\rho$ with no
further factor; the optional factor $(\inf_\vartheta\sym{s}_\mu/\bar{\sym{x}})^{-1}\le1$ of
Lemma~\ref{lem:toeplitz} tightens the bound at the cost of one symbol evaluation.  Because
$\sym{s}_\mu-\rf\,\sym{m}$ is proportional to $C$, $\beta_\mu=C\tilde\beta(G,M,Y)$, and since
$I_a(\xi)$ is decreasing in $a$ (the factor $e^{-ax}$), $\tilde\beta$ is nonincreasing in $G$ and in
$M$; the cells are therefore taken in $(C,Y,\min(G,M))$, not in $\rf$, the minimum over $C$ and over
$(G,M)$ is attained at the cell's corner $(C_{\min},G_{\max},M_{\max})$ exactly, and only the
one-dimensional minimum over $Y$ is taken on a grid with the safety factor $0.99$.  The offline/online evaluation of a residual norm is
subject to cancellation when the residual is small relative to its constituents
\cite{CasenaveErnLelievre2014,BuhrEngwerOhlbergerRave2014}; the implementation uses the
$R$-factor form of these references with an explicit rounding term rather than the Gram-matrix form.
\end{remark}

Theorem~\ref{thm:eim-cert} is what allows a reduced model to be built on $A_Q$ while the
certificate refers to the true operator $\Ah$; the price is the allowance
$\varepsilon_Q\alpha_{\rm lb}^{-1/2}\norm{\overline{\bar\cc}_\theta}_2$ per step, which the choice of $Q$ below
keeps under a tenth of the target tolerance.
\subsection{Candidates, initial target and training objective}
\label{ssec:target}\label{ssec:architectures}\label{ssec:exact-ic}\label{ssec:losses}\label{sec:data}
A candidate predicts coefficients in the fixed spline space, $(\mu,T,\psi,\tau)\mapsto
\cc_\theta(\tau;\mu,T,\psi)$ and $U_{\theta,h}=\sum_k c_{\theta,k}\phi_{h,k}$; its prescribed coefficients are overwritten at every time level before certification, which makes
it admissible in the same affine space as the reference, and its residuals use the assembled
operator of Section~\ref{ssec:assembly}.  The initial vector is the solver's specified target: the
raw mass projection $\bar M\bar\cc^0=\bar\bb^0$ with $\bar b_k^0=(\psi,\phi_{h,k})_{L^2(\Omega)}$,
followed by the imposition of the prescribed coefficients, with free part $\zz_h^0=E^\top\cc_h^0$; all
reference and training trajectories use it, and the initial discrepancy
$\beta_0=\|\zz_h^0-\zz_\theta^0\|_{M_I}$ is evaluated explicitly for every candidate (it vanishes
when the candidate embeds the target, $\cc_\theta(\tau)=\cc_h^0+\tau\mathcal N_\theta$, and is the
projection error for reduced trajectories).  A DeepONet trained on the principal-period
discretisation of the earlier version of this work is certified as a candidate of the exact system
in Section~\ref{sm:deeponet} of the supplementary materials and rejected at every tolerance with no
false accept; the neural candidate of the production study is the POD-NN of
Section~\ref{ssec:pod-construction}.

The training objective combines the relative mass-weighted data loss
$\mathcal L_{\rm data}=(\cc_\theta-\cc_h)^\top\Mh(\cc_\theta-\cc_h)/(\cc_h^\top\Mh\cc_h+\epsilon)$,
the relative squared $L^2(\Omega)$ error of the homogeneous spline field with a small
$\epsilon>0$, and the residual regulariser
\begin{equation}
 \mathcal L_{\rm res}
 =\frac{1}{N_T}\sum_{n=0}^{N_T-1}
 (\rr_\theta^{n+1/2})^\top\Sh(\mu)^{-1}\rr_\theta^{n+1/2},
 \qquad \mathcal L=\mathcal L_{\rm data}+\lambda_R\mathcal L_{\rm res},
 \label{eq:residual-loss}
\end{equation}
with $\rr_\theta^{n+1/2}$ given by \eqref{eq:cn-residual}; on a uniform grid
$\etaSharp^2=T\mathcal L_{\rm res}$, so the regulariser controls the absolute mass error together
with the initial term, whereas the data loss is normalised.  The training set consists of 800
Sobol-sampled trajectories (seed 0) of the exact-assembly reference at the production resolution,
with a separate 30-query validation set (seed 999999) used for checkpoint selection and a
4000-trajectory set for the data-size study; Table~\ref{tab:training} in the supplementary materials
lists the optimiser, schedule, parameter counts and wall times.

\subsection{Mass-weighted POD trajectories}
\label{ssec:pod-construction}
With $M_I=LL^\top$ and the uncentred free snapshot matrix $Z$, define
$V=L^{-\top}U_r$, where $U_r$ contains the leading left singular vectors of
$L^\top Z$. Then $V^\top M_IV=I$ and the subspace minimises the aggregate
training-snapshot projection error in the mass norm. The 800 trajectories
supply 192,800 snapshots ($241$ time levels each). The reduced equations are
\begin{equation}
 \frac{\alpha^{n+1}-\alpha^n}{\dt}
 + (V^\top A_I V)\frac{\alpha^{n+1}+\alpha^n}{2}
 = V^\top f_I^{n+1/2},\qquad
 \alpha^0=V^\top M_I \zz_h^0.
 \label{eq:repaired-pod-cn}
\end{equation}

The reconstructed free candidate $V\alpha^n$ shares the prescribed boundary
coefficients. Its finite-rank initial discrepancy $\beta_0=\norm{(I-VV^\top M_I)\zz_h^0}_{M_I}$
is the projection error of the initial target and enters the same certificate as the neural
candidate.

\paragraph{POD--Galerkin with the affine operator.}
With $A_Q$ in place of $A_I$ the reduced operator is
\[
V^\top A_Q(\mu)V=-(\rf-\omega_S)\,V^\top D_IV+\lambda_0 I-C\Gamma(-Y)h^{-Y}\sum_q\ell_q(\nu)\,V^\top T_qV,
\]
a combination of $Q+2$ precomputed $r\times r$ matrices with scalar coefficients, and the
reduced forcing is affine in the same scalars and in the two lifting values $b_L(\tau_n)$, through
the precomputed vectors $V^\top M_I e_{\rm bd}$, $V^\top T_q e_{\rm bd}$.  The online solve therefore
costs $O(Qr^2+N_Tr^2)$ operations and never touches the full space.  The candidate is certified
either at full order, by Theorem~\ref{thm:eim-cert} with $\eta^X$ or $\eta^S$ and the allowance
term, or online as follows.

\paragraph{Reduced online certificate.}
Fix a comparison matrix $\bar X=r_{\min}M_I+\bar\beta\,h^2R^{(1)}_I$ per cell of a coarse
partition of the $(C,Y,\min(G,M))$ range, with $\bar\beta=0.99\min_{\rm cell}\beta_\mu$ evaluated on
a grid of the cell (Remark~\ref{rem:reduced-online}).  Then $\rho^\top\bar X^{-1}\rho$ is a quadratic form in the
$Q+2$ scalar coefficients and the reduced states, whose $(Q+2)^2$ Gram blocks
$R_q^\top\bar X^{-1}R_{q'}$, $R_q\in\{M_IV,\,T_qV,\,D_IV,\,\text{forcing directions}\}$, are
computed offline by banded solves.  Online, $\eta^{\bar X}$ costs $O(\varrho\,(Q+2)(2r+1))$ operations
with $\varrho$ retained rows (below), independently of $N_h$, and is an upper bound for the exact dual norm without any further factor
(Remark~\ref{rem:reduced-online}); the symbol quantities that accompany it, $\varepsilon_Q(\mu)$ and
the optional tightening factor, are evaluated on an $N$-independent graded grid with derivative
margins.  Two implementation facts matter.  First, the Gram-tensor form of the offline/online
decomposition is numerically unusable here: for accurate candidates the residual is $10^{-8}$ of
the terms it is assembled from, and the Gram form loses all digits.  The implementation uses the
$R$-factor form \cite{CasenaveErnLelievre2014,BuhrEngwerOhlbergerRave2014}, compressed offline by a
singular value decomposition with a rigorous term for the discarded part; the online cost is then
$O(\varrho\,(Q+2)(2r+1))$ with $\varrho$ retained rows ($185$--$1022$ for $r=8$--$128$ at $N=1024$;
$\varrho$ saturates below $N_h$, which is what makes the cost independent of $N$), the online
value agrees with the full-order $\rho^\top\bar X^{-1}\rho$ to $2\times10^{-10}$ for POD--Galerkin
candidates (the double-precision floor, confirmed in extended precision), and the compressed bound
exceeds the uncompressed one by at most $5\times10^{-4}$ relative.  Second, the cell norm costs a
factor $1.87$ (median) against the query's own $X_I$ and about $4$ against $\eta^S$
(Table~\ref{tab:pod-ranks}).  On an idle core the complete online certificate takes
$2.9$, $3.9$, $4.6$, $5.6$, $6.8$, $8.3$~ms at $r=8$, $48$, $64$, $80$, $96$, $128$ ($Q=36$), of which
$2.7$~ms are rank-independent; the reduced solve takes $0.7$--$2.0$~ms, and the offline tensors
occupy $9$~MB ($r=8$) to $688$~MB ($r=128$) for the 18 cells.

\paragraph{POD-NN.}
The neural reduced candidate maps $(\mu,T)$ to the $r\times N_T$ matrix of POD coefficients
$\alpha^1,\dots,\alpha^{N_T}$ by a multilayer perceptron, with $\alpha^0=V^\top M_I\zz_h^0$ fixed
so that $\beta_0$ is again the projection error only
\cite{HesthavenUbbiali2018}.  Its training objective is the least-squares loss in the reduced
coordinates (the mass-weighted loss $\mathcal L_{\rm data}$, since $V^\top M_IV=I$) plus
$\lambda_R$ times the reduced online certificate $\eta^{\bar X}_{\rm aug}$ of the previous paragraph
(not squared), which is a differentiable function of the network output and costs nothing per
sample; it plays the role of the residual regulariser \eqref{eq:residual-loss}, with the exact dual
norm replaced by a rigorous majorant.  The network is evaluated in double precision: single-precision rounding of the
predicted coefficients alone produces residuals whose certified bound is of the order $10^{-3}$ at
$N=512$, i.e.\ the tolerance itself.  The standard output has $r\times N_T=11\,520$ coefficients ($3.09$~M parameters, of which $2.96$~M in
the output layer).  A second variant predicts the coefficients of a temporal POD of the reduced
trajectories, $r\times r_t$ with $r_t=23$ modes chosen so that the temporal projection error is
below $\varepsilon_N/10$ ($0.42$~M parameters); it removes the time-difference noise that the
residual penalises by construction.  Both are trained with $\lambda_R\in\{0,10^{-2},10^{-1},1\}$ on
800 and 4000 trajectories (Table~\ref{tab:training}).

\subsection{Cost model}
\label{ssec:cost}
Let ${\rm it}_B$ be the mean number of preconditioned GMRES iterations per Crank--Nicolson step
(Table~\ref{tab:pcg-v9}: $5$--$9$), and count one FFT of length $2N_h$ as $5N_h\log_2N_h$
flops.  Per query,
\begin{align*}
 \text{direct solve:}&\quad N_T\,{\rm it}_B\,(2\cdot5N_h\log_2N_h+O(N_h)),\\
 \text{certificate }(\eta^X):&\quad N_T\,(2\cdot5N_h\log_2N_h)+N_T(2p-1)N_h+O(p^2N_h),\\
 \text{certificate }(\eta^S):&\quad N_T\,{\rm it}_{\rm PCG}\,(2\cdot5N_h\log_2N_h)+\text{the above},\\
 \text{reduced solve:}&\quad (Q+2)r^2+N_T\,r^2,\\
 \text{reduced certificate }(\eta^{\bar X}):&\quad \varrho\,(Q+2)(2r+1)+N_T\,r\,(Q+2)\quad(\varrho\le N_h\ \text{retained rows}),
\end{align*}
plus the symbol evaluation $O(LK)$ common to all certified variants.  Two conclusions are
independent of the implementation.  First, at full order the certificate $\eta^X$ costs one
residual evaluation per step, about $1/{\rm it}_B$ of the direct solve: in one space dimension
with FFT-applied Toeplitz operators the certificate can be at most $5$--$10$ times cheaper than
the solve it certifies, whatever the candidate.  Second, the reduced solve and certificate do not
involve $N_h$ once the offline compression has fixed $\varrho$: a certified reduced query costs
$O(\varrho\,Q\,r)$ operations at any resolution, with the rigorous link to the full-order problem
supplied by Theorem~\ref{thm:eim-cert} and Lemma~\ref{lem:toeplitz}.  Table~\ref{tab:cost} of the
supplementary materials gives the measured constants.

\section{Numerical evaluation}
\label{sec:numerics}
\subsection{Experimental design}
\label{ssec:design}\label{tab:populations}
All studies use the exact assembly of Section~\ref{ssec:assembly} ($K=100$), quadratic splines,
the box of \eqref{eq:numerical-box} and the free system.  The verification studies use 60
Crank--Nicolson steps at several $N$ on 32 box vertices and 256 Sobol points (288 queries) and on 48
audit queries of the earlier version with perturbed initial vectors; the production studies from
Section~\ref{ssec:calibration} on use the configuration fixed there, $N=1024$, $N_T=240$,
$\varepsilon_N=1.18\times10^{-3}$, and report every verdict at the sweep
$\varepsilon_N\times\{0.1,0.3,1,3,10\}$.  Their populations are the 30 price-diagnostic queries and
two box corners of the calibration (five spots $0.8\mathcal K,\dots,1.2\mathcal K$), the 800 and 4000
Sobol training trajectories (seed 0) with a 30-query validation set (seed 999999) on which the
reduced models and candidates are studied, 120 Sobol benchmark queries, $10^4$ Sobol queries for the
multi-query study, and 200 Sobol queries outside the box plus 64 extended vertices for the
out-of-distribution study; statistics refer to the indicated study only.  The implementation is the
package \texttt{cgmy\_cert} (Section~\ref{ssec:archive}); every number below is read from its result
files.  Statistics over spots first take the maximum within each query, and bootstrap intervals
resample queries, rather than spots or repeated timings.
\subsection{Verification of the assembly, the bounds and the solvers}
\label{ssec:verification}\label{ssec:e0-assembly}\label{ssec:e0-spectral}\label{ssec:e0-riesz}\label{ssec:e0-initial}
The verification studies are reported in full in Section~\ref{sm:verification} of the supplementary
materials; four facts are used in the sequel.  The assembled free matrix agrees with the exact
Galerkin matrix, computed independently at 40 digits, to $10^{-8}$ in operator norm at $K=100$, the
truncation error decreases like $K^{-(2p-Y-1)}$ as Lemma~\ref{lem:tail} predicts, and the bound
$\varepsilon_K$ is within a factor $1.1$--$3.7$ of the observed error at $K=100$
(Table~\ref{tab:assembly-v9}).  The bound of Theorem~\ref{thm:spectral} and the two-sided comparison
\eqref{eq:comparison} hold at all $289\times3$ tested (query, $N$) pairs, with
$\lambda_{\min}(S_I)/\smin(\mu)\in[1.003,17.8]$; the a-priori factor $\sqrt{\Thmu}$ grows like
$(\min(G,M)h)^{(Y-2)/2}$ (median $2.05$ at $N=1024$), while the observed ratio $\eta^X/\eta^S$ on
candidates stays below two at the production resolution (Figures~\ref{fig:spectral-bound}
and~\ref{fig:theta}).  Circulant-preconditioned conjugate gradients solve the Riesz systems in
$8$--$12$ iterations for $N=64$--$4096$ against $42$--$1915$ unpreconditioned, and the GMRES solve of
$B_+$ needs $5$--$9$ iterations per step on average (Table~\ref{tab:pcg-v9}).  On exact dynamics with
perturbed initial vectors both augmented bounds hold with effectivity $1.00\pm10^{-6}$, and for
$Y\in\{0.5,0.8,0.95\}$ the bound $\alpha_{\rm lb}=\rf$ is positive and the certificate holds with
effectivity one, so the discrete framework covers the finite-variation regime without change.
\subsection{Resolution and tolerance calibration}
\label{ssec:calibration}
The production resolution is the smallest $N$ at which the discrete put price is within one cent
of the COS reference \cite{FangOosterlee2008} (8192 modes, range parameter 30; agreement with
Carr--Madan to $2.6\times10^{-12}$ on the checked rows) at the five spots $0.8\mathcal K,\dots,1.2\mathcal K$, over the
30 price-diagnostic queries and the two short-maturity, weakly diffusive box corners
($C=0.05$, $Y=1.1$, $T=0.1$) that had price gaps above two dollars at $N=64$.
Table~\ref{tab:budget} gives the maxima.  The 30 queries alone would pass at $N=512$ ($0.11$ cent),
the corners do not ($2.4$ and $3.2$ cent); at $N=1024$ all are at most $0.235$ cent with 240 steps,
whereas 60 steps give $4.6$ cent on the 30 queries at that resolution.  With the spot factor
$\max_s\beta_I(x_s)=8.44$ of \eqref{eq:gamma-infty}, one cent corresponds to the trajectory
tolerance $\varepsilon_N=1.18\times10^{-3}$; this configuration is fixed once and not changed by
any later experiment.  A tenth of a cent ($\varepsilon_N/10$) is reported as a certificate
tolerance only, since the discretisation error at $N=1024$ is $0.235$ cent.
\begin{table}[!htbp]
\centering\small
\caption{Resolution and tolerance calibration: maximum price error $\max|V_h(T,S)-V_{\mathrm{COS}}|$ over the 30 price-diagnostic queries and the two box corners, spots $S\in\{0.8,\dots,1.2\}\mathcal K$; COS with 8192 modes, range parameter 30. $\varepsilon_N=\tau_\$/\max_s\beta_I(x_s)$ with $\tau_\$=1$ cent.}
\label{tab:budget}
\begin{tabular}{rrrrrrr}
\toprule
$N$ & $N_T$ & 30 queries & corner A & corner B & $\max_s\beta_I$ & $\varepsilon_N$\\
\midrule
128 & 60 & $0.0361$ & $0.703$ & $0.783$ & $3.103$ & $3.22\times10^{-3}$\\
128 & 240 & $0.0361$ & $0.703$ & $0.783$ & $3.103$ & $3.22\times10^{-3}$\\
256 & 60 & $2.82\times10^{-3}$ & $0.172$ & $0.203$ & $4.396$ & $2.28\times10^{-3}$\\
256 & 240 & $2.80\times10^{-3}$ & $0.172$ & $0.203$ & $4.396$ & $2.28\times10^{-3}$\\
512 & 60 & $6.17\times10^{-3}$ & $0.0236$ & $0.0316$ & $6.073$ & $1.65\times10^{-3}$\\
512 & 240 & $1.12\times10^{-3}$ & $0.0236$ & $0.0316$ & $6.073$ & $1.65\times10^{-3}$\\
1024 & 60 & $0.0464$ & $1.36\times10^{-3}$ & $2.36\times10^{-3}$ & $8.441$ & $1.18\times10^{-3}$\\
1024 & 240 & $1.59\times10^{-3}$ & $1.35\times10^{-3}$ & $2.35\times10^{-3}$ & $8.441$ & $1.18\times10^{-3}$\\
\bottomrule
\end{tabular}
\end{table}

\subsection{Reduced model and neural candidates}
\label{ssec:e1}\label{ssec:pod-results}
\begin{figure}[!htbp]
\centering
\includegraphics[width=0.78\linewidth]{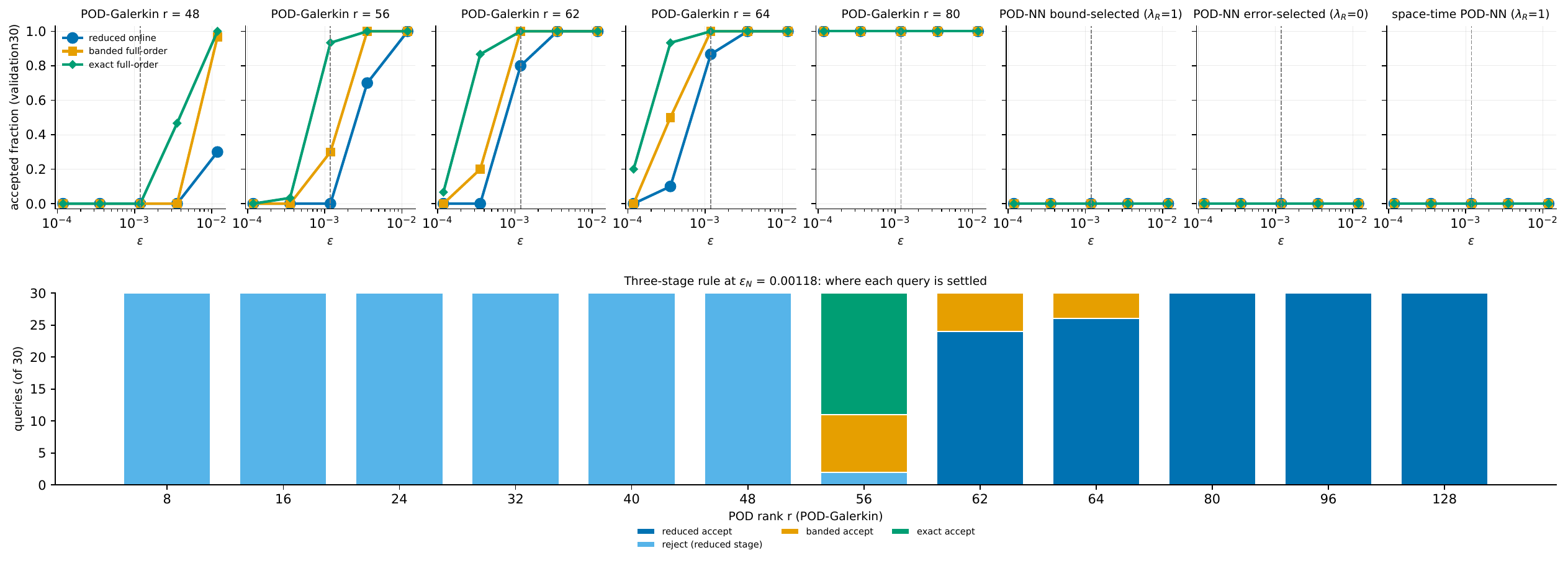}
\caption{Top: accepted fraction of the 30 validation queries as a function of the tolerance
$\varepsilon$ (dashed line: $\varepsilon_N$) for POD--Galerkin candidates of rank $48$--$80$ and the
two best POD-NN checkpoints, under the reduced online, banded full-order and exact full-order
certificates.  Bottom: the stage at which the three-stage rule of Proposition~\ref{prop:threestage}
settles each query at $\varepsilon_N$, per rank.}
\label{fig:acceptance}
\end{figure}

Table~\ref{tab:pod-ranks} and Figure~\ref{fig:acceptance} report the POD--Galerkin candidates on
the validation set at the production configuration.  The mass-weighted basis from 800 trajectories
reaches a median error of $3.1\times10^{-4}$ at $r=48$ and $1.7\times10^{-5}$ at $r=64$; five times
more training data does not enlarge the space noticeably (singular values within $10\%$ up to $r=48$, within a factor $1.5$ up to $r=128$).  The
exact certificate has effectivity $9$--$15$ for $16\le r\le80$ ($6.8$ at $r=8$), consistent with Theorem~\ref{thm:efficiency}
(efficiency holds in the energy norm, which weights the high-frequency part of the error far more
than the mass norm in which the error is measured); the banded certificate is a further factor
$\approx2$ above it, as predicted by $\sqrt{\Thmu}$, and the reduced online certificate a further
factor $\approx2$, the cell factor.  Acceptance at one cent therefore starts at $r=56$ under the
exact certificate (28 of 30), at $r=62$--$64$ under all three, and is complete at $r\ge80$.  The
three-stage rule settles every query at $r\le48$ by cheap rejection and every query at $r\ge80$ by
reduced acceptance; the full-order stages are reached only in the band $r=56$--$64$.  No false
accept occurs under any certificate at any tolerance.  At $r\ge96$ the true error ($4\times10^{-8}$)
falls below the floor set by the EIM allowance at $Q=36$, which raises the reduced effectivity to
$230$--$420$ without affecting acceptance.
\begin{table}[!htbp]
\centering\small
\caption{POD--Galerkin candidates on the 30 validation queries at the production configuration ($N=1024$, $N_T=240$, $Q$ matched to $\varepsilon$): median maximum-in-time mass-norm error against the direct solution, median effectivity (bound/error) of the three certificates, queries accepted at $\varepsilon_N=1.18\times10^{-3}$ under each certificate, and the stage at which the three-stage rule of Proposition~\ref{prop:threestage} settles the query (reduced accept / cheap reject / banded accept / exact accept / exact reject). Basis from 800 training trajectories; the 4000-trajectory basis gives the same picture (Section~\ref{ssec:e1}). No false accepts occur under any certificate at any tolerance of the sweep.}
\label{tab:pod-ranks}
\resizebox{\textwidth}{!}{%
\begin{tabular}{rrrrrrr}
\toprule
$r$ & error & eff.\ $\eta^S$ & eff.\ $\eta^X$ & eff.\ $\eta^{\bar X}$ & accepted $\bar X$/$X$/$S$ & three-stage\\
\midrule
8 & 1.19 & 6.8 & 10.8 & 24 & 0/0/0 & 0/30/0/0/0\\
16 & 0.139 & 14.9 & 25.2 & 50 & 0/0/0 & 0/30/0/0/0\\
24 & 0.0298 & 13.9 & 25.2 & 52 & 0/0/0 & 0/30/0/0/0\\
32 & $7.1\times10^{-3}$ & 13.3 & 23.2 & 48 & 0/0/0 & 0/30/0/0/0\\
40 & $1.3\times10^{-3}$ & 13.2 & 24.1 & 47 & 0/0/0 & 0/30/0/0/0\\
48 & $3.1\times10^{-4}$ & 11.7 & 21.6 & 49 & 0/0/0 & 0/30/0/0/0\\
56 & $6.0\times10^{-5}$ & 11.5 & 21.5 & 47 & 0/9/28 & 0/2/9/19/0\\
62 & $1.9\times10^{-5}$ & 11.3 & 22.0 & 50 & 24/30/30 & 24/0/6/0/0\\
64 & $1.7\times10^{-5}$ & 10.5 & 20.0 & 42 & 26/30/30 & 26/0/4/0/0\\
80 & $6.6\times10^{-7}$ & 9.5 & 19.3 & 47 & 30/30/30 & 30/0/0/0/0\\
96 & $3.6\times10^{-8}$ & 9.1 & 16.8 & 231 & 30/30/30 & 30/0/0/0/0\\
128 & $9.9\times10^{-9}$ & 4.3 & 5.3 & 417 & 30/30/30 & 30/0/0/0/0\\
\bottomrule
\end{tabular}}
\end{table}

The neural candidates are reported in Table~\ref{tab:podnn} of the supplementary materials.  None passes at $\varepsilon_N$, nor
at $10^{-3}$: the best standard-output network ($4000$ trajectories, $\lambda_R=1$) has a median true
error of $0.46$ and a median exact bound of $1.26$, about $10^3\varepsilon_N$; the space--time
variant reduces both by an order of magnitude (errors $0.014$--$0.016$ and bounds $0.054$--$0.058$ over three seeds,
$46$--$49\,\varepsilon_N$).  Three observations explain the result.  With $\lambda_R=0$ the network
fits the coefficients best (error $0.08$) but its trajectory is noisy in time, and the residual,
which contains the time difference divided by $\dt$, gives a bound of $12$ (effectivity $150$);
with $\lambda_R>0$ the bound drops by a factor ten while the error rises by a factor five, so the
residual term buys temporal smoothness with accuracy, and selecting checkpoints by the rigorous
bound picks the smoother, less accurate network (bound $1.26$, error $0.46$, against bound $12.2$,
error $0.082$ when selecting by the true error, which the certified user cannot do).  The
space--time representation removes the temporal noise by construction and wins on both measures.
The projection of the reference trajectories onto the same basis, certified at
$3.2\times10^{-3}$, shows that the basis is not the obstacle; the selected epochs ($1890$--$2000$ of
$2000$) show that the networks had not converged, but longer training cannot close a gap of
$46$--$1000$.  Every neural query is settled by the cheap rejection of the first stage, at a cost of
milliseconds and with no false accept.  Figure~\ref{fig:online-time} shows the online cost of the
reduced model on the validation queries: the reduced solve takes $0.7$--$2.0$~ms and the reduced
certificate $2.9$--$8.3$~ms for $r=8$--$128$; the three-stage total at $\varepsilon_N$ equals the
reduced certificate wherever the first stage decides ($r\le48$ by cheap rejection, $r\ge62$ by
acceptance) and rises to $154$~ms only at $r=56$, where the full-order stages are needed.
\subsection{Matched benchmark and cost under refinement}
\label{ssec:benchmark}\label{ssec:scaling-results}
\begin{figure}[!t]
\centering
\includegraphics[width=0.76\linewidth]{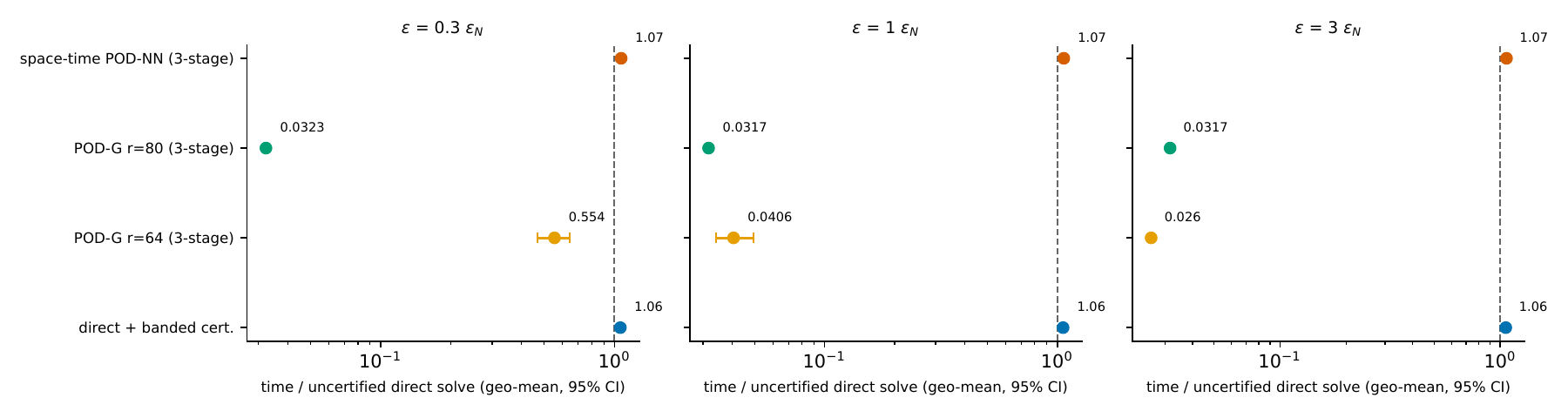}\\[2mm]
\includegraphics[width=0.76\linewidth]{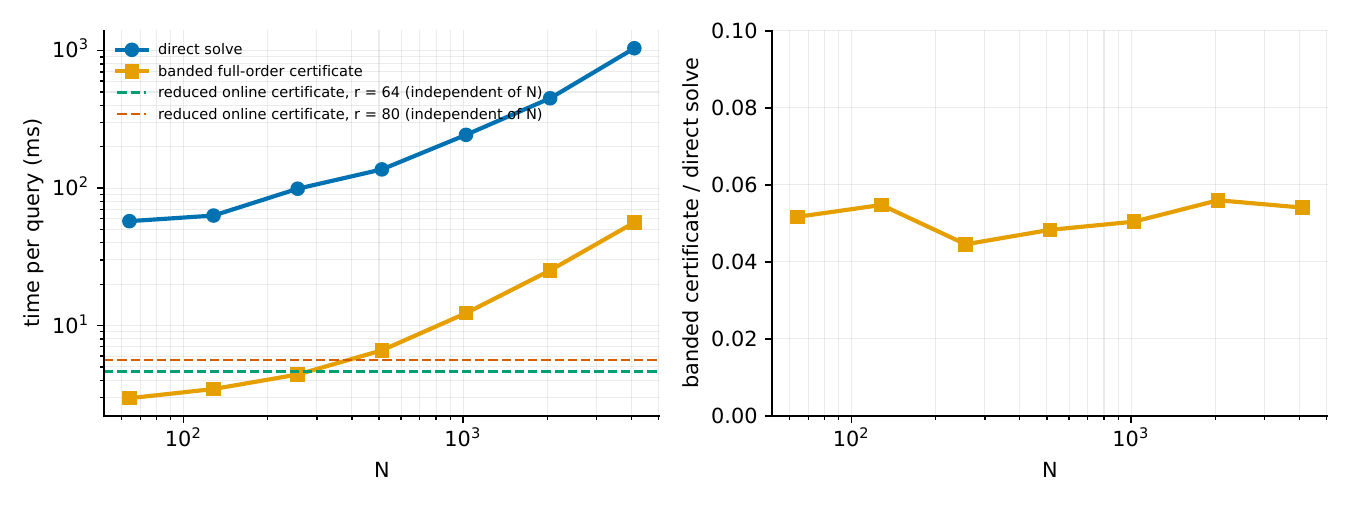}
\caption{Top: time ratio to the uncertified direct solve on the 120 benchmark queries at
$\varepsilon=0.3\varepsilon_N$, $\varepsilon_N$, $3\varepsilon_N$ (geometric means, $95\%$
query-bootstrap intervals, reusable setup).  Bottom left: direct solve and banded full-order
certificate per query as functions of $N$ (20 benchmark queries, $N_T=240$), with the reduced
online certificate at $r=64$ and $80$, which does not depend on $N$.  Bottom right: the certificate
as a fraction of the solve, $4.5$--$5.6\%$ for $N=64$--$4096$.}
\label{fig:benchmark}
\end{figure}

Table~\ref{tab:benchmark} and Figure~\ref{fig:benchmark} report the 120 benchmark queries at the
production configuration.  The baseline is the \emph{uncertified} circulant-preconditioned GMRES
solve of Section~\ref{ssec:e0-riesz} ($259$~ms with reusable setup, $337$~ms cold); every ratio is a
geometric mean over queries with a $95\%$ query-bootstrap interval, and the certificate cost of the
direct solve is shown separately.  Three results.

First, the full-order certificate is cheap: the banded certificate $\eta^X$ with the exact
operator costs $6.0\%$ of the solve ($7.2\%$ including the per-query spectral data in the fast
online mode of Section~\ref{ssec:spectral}), so a certified direct solve takes $1.06$ times an
uncertified one, and this fraction stays at $4.5$--$5.6\%$ from $N=64$ to $N=4096$
(Figure~\ref{fig:benchmark}, right), as the cost model of Section~\ref{ssec:cost} predicts
(Table~\ref{tab:cost}).  The spectral quantities evaluated on the fine verification grid of
Section~\ref{ssec:e0-spectral} would instead cost $478$~ms and make the certified solve $2.92$
times the uncertified one; the fast mode (a coarse graded grid of $1012$ nodes with ten aliases and
a rigorous $\sin^{2p}$-weighted tail) reproduces $\beta_\mu$ and $\Thmu$ to $0.2\%$ and changes
$\eta^X$ by $0.002\%$.

Second, the certified reduced model is an order of magnitude faster than the direct solve
whenever the three-stage rule stops at the reduced stage: POD--Galerkin with $r=80$ takes
$0.032$ of the direct time ($8.3$~ms) at every tolerance of the sweep, with all 120 queries accepted
at the $N$-independent stage; with $r=64$ it takes $0.041$ at $\varepsilon_N$ (102 reduced and 18
banded accepts), $0.026$ at $3\varepsilon_N$, and $0.55$ at $0.3\varepsilon_N$, where most queries
need a full-order stage and four need correction (without reusable setup the ratio rises to $0.30$
because the query-independent initial projection, $90$~ms, is rebuilt inside the timer).
Third, the neural candidate costs $1.07$ times the direct solve at every tolerance: every query is
rejected at the reduced stage in a few milliseconds and then corrected by a full-order march
warm-started from the candidate and re-certified.  No false accept occurs for any method or tolerance.
\begin{table}[!htbp]
\centering\small
\caption{Matched benchmark on the 120 benchmark queries at the production configuration (one thread, idle machine, reusable setup, five repetitions): geometric-mean time ratio to the uncertified direct solve with 95\% query-bootstrap interval, three-stage outcome (reduced accept / cheap reject / banded or exact accept / corrected), and false accepts.  The direct solve takes $259$~ms; the certificate proper is $6.0\%$ of it and $7.2\%$ with the per-query spectral data.}
\label{tab:benchmark}
\begin{tabular}{llrcr}
\toprule
$\varepsilon$ & method & time ratio (CI) & outcome & f.a.\\
\midrule
-- & direct $+$ banded certificate & 1.06 (1.06--1.06) & cert.\ fraction 0.07 & --\\
$0.3\varepsilon_N$ & POD--Galerkin $r=64$ & 0.554 (0.47--0.642) & 8/4/108/4 & 0\\
$0.3\varepsilon_N$ & POD--Galerkin $r=80$ & 0.0323 (0.032--0.0326) & 120/0/0/0 & 0\\
$0.3\varepsilon_N$ & space--time POD-NN & 1.07 (1.06--1.07) & 0/120/0/120 & 0\\
$1\varepsilon_N$ & POD--Galerkin $r=64$ & 0.0406 (0.034--0.0494) & 102/0/18/0 & 0\\
$1\varepsilon_N$ & POD--Galerkin $r=80$ & 0.0317 (0.0314--0.032) & 120/0/0/0 & 0\\
$1\varepsilon_N$ & space--time POD-NN & 1.07 (1.06--1.07) & 0/120/0/120 & 0\\
$3\varepsilon_N$ & POD--Galerkin $r=64$ & 0.026 (0.0257--0.0262) & 120/0/0/0 & 0\\
$3\varepsilon_N$ & POD--Galerkin $r=80$ & 0.0317 (0.0313--0.032) & 120/0/0/0 & 0\\
$3\varepsilon_N$ & space--time POD-NN & 1.07 (1.06--1.07) & 0/120/0/120 & 0\\
\bottomrule
\end{tabular}
\end{table}

\subsection{Out-of-distribution detection}
\label{ssec:ood}
The certificate makes no assumption about where a query lies, so the natural test is to leave the
box.  Table~\ref{tab:ood} reports 264 queries outside the box \eqref{eq:numerical-box}
(200 Sobol points in the shells $C\in[1,2]$, $Y\in[0.9,1.1)\cup(1.8,1.95]$, $G\in[0.5,1)$,
$M\in[9,15]$, excluding $|Y-1|<0.02$ where $\Gamma(-Y)$ is singular, and the 64 vertices of the
extended box).  Every query is certifiable: $\alpha_{\rm lb}=\rf-\varepsilon_K>0$ and the Toeplitz
aliasing stays below $10^{-15}$ even for $G<1$.  Outside the cells of Section~\ref{ssec:eim} the
three-stage rule starts at the banded stage with the query's own $X_I$, and the EIM error
$\varepsilon_Q$, evaluated at the query, grows by two orders of magnitude (median $2\times10^{-7}$
against $1.7\times10^{-9}$ inside the box) because the interpolant extrapolates in $Y$ and $G$; it
enters the allowance, not the verdict.

\begin{figure}[!htbp]
\centering
\includegraphics[width=0.72\linewidth]{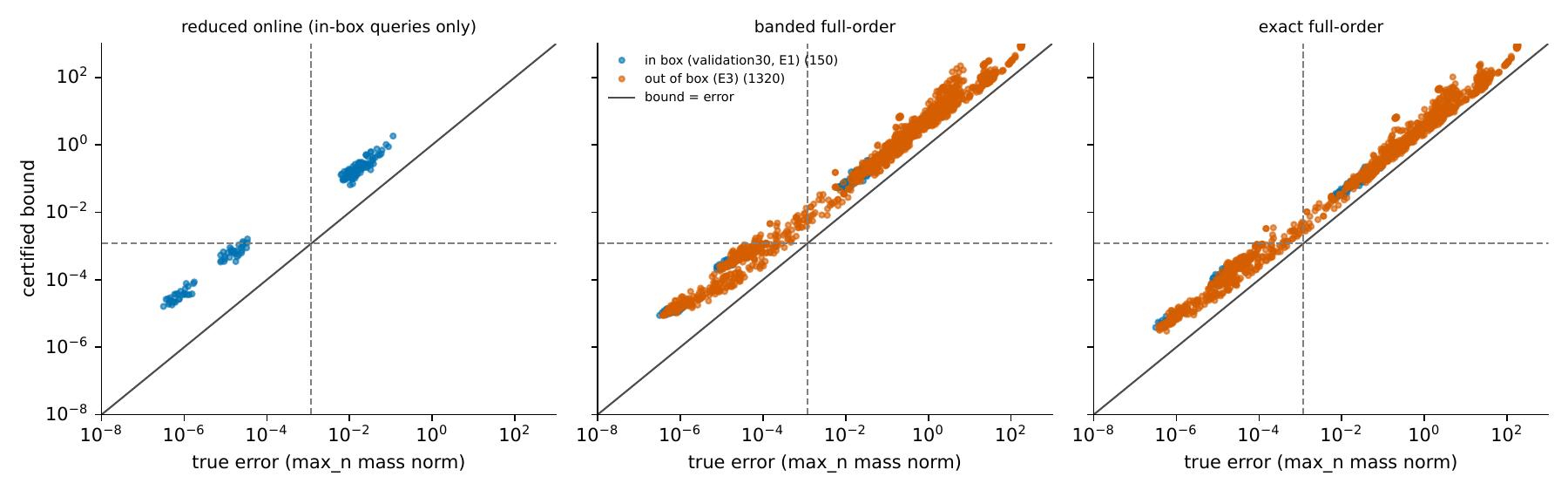}
\caption{Certified bound against true error for the 264 out-of-box queries and the 30 validation
queries (POD--Galerkin $r=64$ and $80$, space--time POD-NN), under the three certificates; the lines
mark $\varepsilon_N$ and bound $=$ error.}
\label{fig:ood}
\end{figure}

Three facts follow.  There is no false accept for any candidate under any certificate at any
tolerance of the sweep, as Theorem~\ref{thm:fully-discrete-certificate} requires.  The reduced
model keeps working on most of the extended domain: POD--Galerkin with $r=80$ is within one cent
on 210 of the 264 queries, and the certificate accepts 193 of them (92\%), with 17 false rejects
from the effectivity of the exact certificate (median $6.3$); its failures concentrate in the shells $Y>1.8$ and $G<1$, which
are also the shells of largest $\varepsilon_Q$.  A detector that rejects everything outside the
training box has the same recall as the certificate but discards all 210 working queries
(precision $0.20$ against $0.76$); the seed-ensemble disagreement of the neural candidate, the usual
heuristic, cannot even be scored here because every space--time POD-NN error exceeds
$\varepsilon_N$ in and out of the box.  In Figure~\ref{fig:ood} every point lies above the diagonal.
\begin{table}[!htbp]
\centering\small
\caption{Out-of-distribution study: 264 queries outside the parameter box (200 Sobol points in the
shells $C\in[1,2]$, $Y\in[0.9,1.1)\cup(1.8,1.95]$, $G\in[0.5,1)$, $M\in[9,15]$, and 64 vertices of the
extended box) at the production configuration.  Top: queries on which the candidate's true error
is below $\varepsilon_N$ (``works''), queries the exact certificate accepts (``certified''), and
precision of two rules as detectors of ``true error $>\varepsilon_N$'': the certificate (recall 1 by
theorem) and the naive rule that rejects every query outside the training box (recall 1, but it
discards every working query).  Bottom: POD--Galerkin $r=80$ by shell (a query leaving the box in
several coordinates counts in each), with the median true error and the median EIM error
$\varepsilon_Q$ (inside the box: $1.7\times10^{-9}$).}
\label{tab:ood}
\begin{tabular}{lrrrr}
\toprule
candidate & works & certified & precision (certificate) & precision (naive rule)\\
\midrule
POD--Galerkin $r=64$ & 196 & 181 & 0.82 & 0.26\\
POD--Galerkin $r=80$ & 210 & 193 & 0.76 & 0.20\\
space--time POD-NN & 0 & 0 & -- & --\\
\bottomrule
\end{tabular}
\par\smallskip
\begin{tabular}{lrrrrr}
\toprule
shell ($r=80$) & queries & works & certified & median error & median $\varepsilon_Q$\\
\midrule
$C>1$ & 122 & 104 & 95 & $3.7\times10^{-6}$ & $5.3\times10^{-8}$\\
$Y<1.1$ & 39 & 36 & 31 & $6.1\times10^{-6}$ & $1.4\times10^{-8}$\\
$Y>1.8$ & 36 & 18 & 15 & $2.8\times10^{-3}$ & $6.1\times10^{-7}$\\
$G<1$ & 27 & 21 & 13 & $3.9\times10^{-4}$ & $3.1\times10^{-6}$\\
$M>9$ & 111 & 101 & 94 & $8.6\times10^{-6}$ & $2.5\times10^{-7}$\\
vertices & 64 & 34 & 31 & $2.0\times10^{-4}$ & $4.3\times10^{-6}$\\
\bottomrule
\end{tabular}
\end{table}

\subsection{Multi-query cost and amortisation}
\label{ssec:multiquery}
\begin{figure}[!htbp]
\centering
\begin{minipage}[t]{0.49\linewidth}\centering
\includegraphics[width=\linewidth]{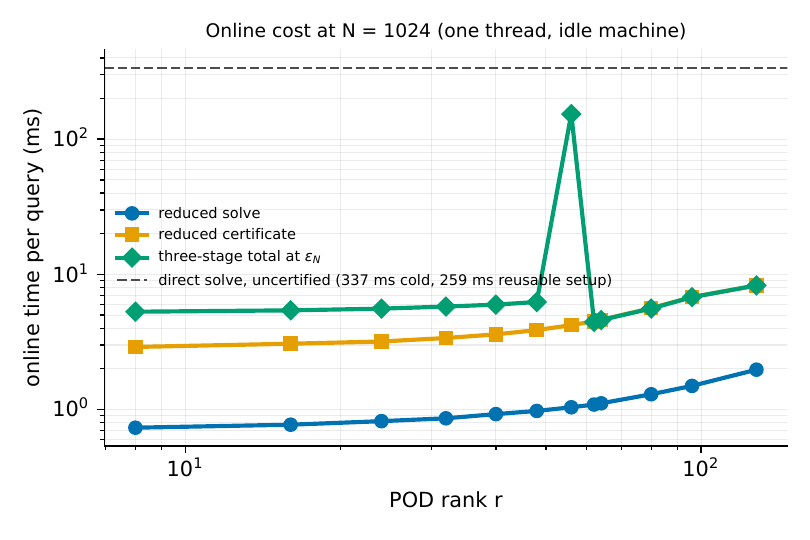}
\caption{Online time per query of the reduced model at the production configuration (one thread,
idle machine): reduced solve, reduced online certificate, and the three-stage total at
$\varepsilon_N$, against the uncertified direct solve (dashed: cold-start and reusable-setup medians of the
benchmark of Section~\ref{ssec:benchmark}).}
\label{fig:online-time}
\end{minipage}\hfill
\begin{minipage}[t]{0.49\linewidth}\centering
\includegraphics[width=\linewidth]{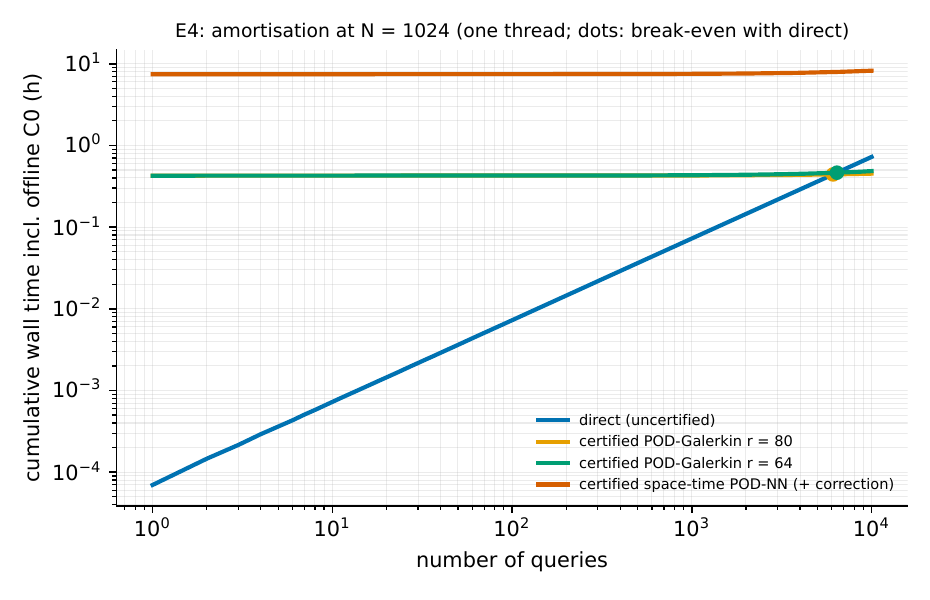}
\caption{Cumulative wall time, offline cost included, against the number of queries at the
production configuration (one thread); dots mark the break-even with the uncertified direct
solver.}
\label{fig:amortisation}
\end{minipage}
\end{figure}

Figure~\ref{fig:amortisation} reports $10^4$ Sobol queries at $\varepsilon_N$ on one core, with
every method's offline cost as the intercept.  The direct solve averages $260$~ms per query.
Certified POD--Galerkin with $r=80$ averages $9.0$~ms: $99.77\%$ of the queries are settled at the
reduced stage, $0.14\%$ at the banded stage, $0.02\%$ at the exact stage, and $0.07\%$ are rejected
and corrected; with $r=64$ the mean is $21.2$~ms (median $7.2$~ms; $89.7\%$ reduced, $9.8\%$ banded,
$0.35\%$ corrected).  The offline cost of the reduced model is about $1530$~s: reference
trajectories $208$~s, POD $797$~s, EIM $65$~s, cell constants $425$~s, offline tensors $39$~s
(single-thread equivalents; the POD used four BLAS threads).  The fitted totals
$T(Q)=1533\,\mathrm s+9.0\,\mathrm{ms}\cdot Q$ ($r=80$) and $1533\,\mathrm s+21.2\,\mathrm{ms}\cdot Q$
($r=64$) break even with $260\,\mathrm{ms}\cdot Q$ after about $6100$ and $6400$ queries.  The
space--time POD-NN never breaks even: its online cost exceeds the direct solve ($282$~ms, since
every query is corrected), and its offline cost ($2.7\times10^4$~s) is dominated by training.
No reference solutions exist for these $10^4$ queries, so this experiment reports verdicts and
times only; the false-accept checks are those of Sections~\ref{ssec:e1}--\ref{ssec:ood}.
\subsection{Reproducibility}
\label{ssec:archive}
All computations use the package \texttt{cgmy\_cert} (Python 3.10, NumPy 1.26, SciPy 1.15, mpmath
1.3, PyTorch 2.14 on CPU) on one core of an Intel Core Ultra 9 285K, with the thread count pinned to
one for every timed quantity, a warm-up query, five repetitions in shuffled order and medians of
per-query medians.  The archive contains the package with its test suite (67 tests, including the
40-digit assembly fixture and the 459-coordinate regression against the earlier implementation),
the scripts that regenerate every result file, the result files themselves with hardware, software
versions and commit hashes, and the SHA-256 manifest of the files too large to archive (reference
trajectories, bases, offline tensors, network checkpoints; 27~GB) with the command that regenerates
each.  A clean clone rebuilds every table byte-identically and every figure pixel-identically from
the result files and reruns the verification studies in 32 minutes on one core; every number in
this section is read from a result file by a checked script.
\FloatBarrier

\section{Interpretation and conclusion}
\label{sec:conclusion}\label{sec:discussion}
The discrete problem of this paper is the exact spline Galerkin system of the CGMY operator, and
the first point is that for shift-invariant bases ``exact'' is cheap: the matrices are Toeplitz
sections of a periodised symbol whose alias series is summed to a stated tail, and the same symbol
gives, in closed form, the spectral lower bound, the banded comparison norm and the affine
structure on which everything else rests.  Positivity of the symmetric part, which the earlier
version of this work assessed by an audit of 459 parameter points and a numerical eigenvalue
estimate, is a theorem with the explicit constant $\alpha_{\rm lb}=\rf-\varepsilon_K$, valid for
$Y\in(0,2)\setminus\{1\}$.

The second point is the cost of a guarantee.  At full order the certificate is one residual
evaluation per time step plus a banded solve, about $5\%$ of the direct solve at every resolution
tested; in one space dimension with FFT-applied operators it cannot be much less, because the
solve itself is a handful of residual evaluations per step.  The scalable statement is the reduced
one: with the affine operator and a fixed comparison norm per cell of parameter space, the
certificate of a reduced trajectory costs a few milliseconds independently of $N$, the three-stage
rule escalates to the full-order bounds only for candidates within a known factor of the
tolerance, and on $10^4$ queries $99.8\%$ of the verdicts are settled without touching the full
space, at $0.03$ of the direct time from rank $80$ on.

The third point is the negative one.  No neural reduced model trained here comes within a factor
$46$ of the tolerance, although the basis it predicts in certifies at $3\times10^{-3}$ and the
projection-based reduced model on the same basis passes at rank $56$--$64$.  The certificate
establishes this query by query, with no false accept on any of the $150$ in-box and $264$
out-of-box queries with a reference solution, and explains it: the residual penalises temporal
roughness that a least-squares fit of the coefficients does not see, so the networks that the
rigorous bound selects are the smoother, less accurate ones.  Whether a better architecture, more
data or a different loss closes the remaining gap is a question the framework can now answer at the
price of a few milliseconds per candidate.  Outside the training box the same certificate accepts
$92\%$ ($193$ of $210$) of the queries on which the reduced model still delivers one cent, where
any detector based on the box itself would discard them all.

What the certificate does not do is bound the continuum price: the resolution was chosen so that
the discretisation error is below one cent on the diagnostic queries, and the certified bound is
added to that budget, not instead of it; the grid evaluation of the symbol extrema is not a
validated enclosure (Section~\ref{ssec:spectral}), and the floating-point residual is covered by
an allowance rather than by interval arithmetic (Remark~\ref{rem:arithmetic-scope}).  The American
put is the natural next target: Theorem~\ref{thm:vi} in Section~\ref{ssec:vi} of the supplementary
materials gives the residual certificate for the backward-Euler discrete variational inequality in
the same spline space, with the complementarity term as the only addition, so that a certified
reduced model for the American price is a matter of the obstacle solver and training data, not of
new theory.
\FloatBarrier

\end{document}


\maketitle
\section{Explicit constants and proof of uniform stability}
\label{app:proof-uniform}\label{ssec:lk-coercivity}
For the box \eqref{eq:parameter-set}, define
\begin{equation}
  a_\ast:=\max\{M_{\max},G_{\max}\},\quad
  \kappa_Y:=\frac{2e^{-\pi}}{(2-Y)\pi^{Y}},\quad
  \Theta_w:=\max\{1,M_{\max}^2,G_{\max}^2\}^{Y_{\max}/2}.
  \label{eq:uniform-constants-aux}
\end{equation}
Explicit admissible choices are
\begin{equation}
  c_1
  =\min_{\mu\in\Pset}\min\Bigl\{
    \frac{\rf}{1+2^{1+Y/2}a_\ast^{Y}},\;
    \frac{C\,\kappa_Y}{2^{Y/2}}
  \Bigr\}
  =\frac{r_{\min}}{1+2^{1+Y_{\max}/2}a_\ast^{Y_{\max}}},
  \label{eq:c1-explicit}
\end{equation}
the second branch being larger on the box of equation~(\ref{eq:numerical-box}) of the main text, and
\begin{equation}
  C_{\rm cont}
  =2\Bigl(
    \Omega_\ast
    +2C_{\max}\gamma_{\max}\max\{1,a_\ast^2\}^{Y_{\max}/2}
    +\lambda_\ast
  \Bigr),
  \label{eq:Ccont-explicit}
\end{equation}
where
\begin{equation}
\begin{gathered}
\gamma_{\max}:=\!\!\max_{Y\in[Y_{\min},Y_{\max}]}\!\!\Gamma(-Y)<\infty,
\qquad
\Omega_\ast:=\sup_{\mu\in\Pset}|\rf-\omega_S(\mu)|,
\qquad
\lambda_\ast:=\sup_{\mu\in\Pset}\lambda_0(\mu),
\end{gathered}
\label{eq:uniform-aux-constants}
\end{equation}
the last two being finite.
Throughout, $\mu=(C,G,M,Y,\rf)\in\Pset$ and $s_\mu=Y/2\in(1/2,1)$.  We work on
the full line and transfer to $\Omega$ at the end using
the zero-extension convention
\eqref{eq:energy-space}, which preserves all constants.
\medskip\noindent
\textbf{Step 1: uniform equivalence of the parameter-dependent norm.}
Let $s\in[Y_{\min}/2,Y_{\max}/2]\subset(1/2,1)$ and let
$w_\mu(\xi):=1+|M-i\xi|^{2s}+|G+i\xi|^{2s}
=1+(M^2+\xi^2)^{s}+(G^2+\xi^2)^{s}$ be the weight in
\eqref{eq:HMG-norm}.  Since $s\le1$ and $(a+b)^s\le a^s+b^s$ for $s\in(0,1]$,
\[
  w_\mu(\xi)\ \le\ 1+2\max\{1,M_{\max}^2,G_{\max}^2\}^{s}(1+\xi^2)^{s}
  \ \le\ (1+2\Theta_w)(1+\xi^2)^{s},
\]
with $\Theta_w=\max\{1,M_{\max}^2,G_{\max}^2\}^{Y_{\max}/2}$, using
$(1+\xi^2)^s\ge1$.  For the lower bound, $w_\mu(\xi)\ge1+|\xi|^{2s}$ and,
since $(1+\xi^2)^s\le2^s\max\{1,|\xi|^{2s}\}\le2(1+|\xi|^{2s})$ for
$s\le1$,
\[
  w_\mu(\xi)\ \ge\ \tfrac12(1+\xi^2)^{s}.
\]
Hence, uniformly in $\mu\in\Pset$,
\begin{equation}
  \tfrac12\norm{u}_{H^{s_\mu}(\R)}^2
  \ \le\ \norm{u}_{H^{s_\mu}_{M,G}(\R)}^2
  \ \le\ (1+2\Theta_w)\norm{u}_{H^{s_\mu}(\R)}^2 ,
  \label{eq:uniform-norm-equivalence}
\end{equation}
where $\norm{u}_{H^{s}(\R)}^2=\tfrac1{2\pi}\int(1+\xi^2)^{s}|\widehat u|^2$.

\medskip\noindent
\textbf{Step 2: uniform continuity.}
From \eqref{eq:symbol}, for every $\xi$,
\[
  |A_\mu(\xi)|
  \le|\rf-\omega_S||\xi|
   +C\Gamma(-Y)\bigl[(M^2+\xi^2)^{Y/2}+(G^2+\xi^2)^{Y/2}\bigr]
   +\lambda_0 .
\]
The map $(C,G,M,Y)\mapsto\omega_S$ of \eqref{eq:martingale-correction} is
continuous on $\Pset$ because $(M-1)^Y$ is continuous and finite on the selected box
(and remains so at $M=1$), and $\Gamma(-Y)$ is continuous on $[Y_{\min},Y_{\max}]\subset(1,2)$;
hence $\Omega_\ast<\infty$ and $\lambda_\ast<\infty$.  Since $Y>1$ we have
$|\xi|\le(1+\xi^2)^{Y/2}$, and
$(a^2+\xi^2)^{Y/2}\le\max\{1,a^2\}^{Y/2}(1+\xi^2)^{Y/2}$.  Therefore
\begin{equation}
  |A_\mu(\xi)|
  \le\Bigl(\Omega_\ast+2C_{\max}\gamma_{\max}\max\{1,a_\ast^2\}^{Y_{\max}/2}
   +\lambda_\ast\Bigr)(1+\xi^2)^{s_\mu}.
  \label{eq:uniform-symbol-upper}
\end{equation}
By Plancherel's theorem and the Cauchy--Schwarz inequality,
\[
  |a_\R(u,w)|
  \le\Bigl(\sup_{\xi\in\R}\frac{|A_\mu(\xi)|}{(1+\xi^2)^{s_\mu}}\Bigr)
  \norm{u}_{H^{s_\mu}(\R)}\norm{w}_{H^{s_\mu}(\R)},
\]
and
\eqref{eq:uniform-norm-equivalence} converts the $H^{s_\mu}$ norms into
$H^{s_\mu}_{M,G}$ norms at the cost of a factor $2$.  This gives
\eqref{eq:Ccont-explicit}.

\medskip\noindent
\textbf{Step 3: the L\'evy--Khintchine form.}
By Lemma~\ref{lem:LK},
\begin{equation}
  \Real A_\mu(\xi)=\rf+C\,[I_G(\xi)+I_M(\xi)],\qquad
  I_a(\xi)=\int_0^\infty(1-\cos\xi x)\,e^{-ax}x^{-1-Y}dx\ \ge0,
  \label{eq:uniform-low-freq}
\end{equation}
so $\Real A_\mu(\xi)\ge \rf\ge r_{\min}>0$ for all $\xi$; the drift term is purely
imaginary and does not enter.  (This is where $r_{\min}>0$ is used.)  Two elementary
lower bounds for $I_a$ follow from $1-\cos u\ge 2u^2/\pi^2$ on $|u|\le\pi$: restricting the
integral to $0<x\le\pi/|\xi|$ and substituting $u=ax$,
\begin{equation}
  I_a(\xi)\ \ge\ \frac{2\xi^2}{\pi^2}\,a^{Y-2}\!\int_0^{a\pi/|\xi|}\!u^{1-Y}e^{-u}\,du
  \ \ge\ \frac{2}{\pi^2}\,\gamma_Y\,a^{Y-2}\xi^2\quad(|\xi|\le a),
  \label{eq:LK-low}
\end{equation}
with $\gamma_Y:=\int_0^\pi u^{1-Y}e^{-u}\,du$, and, bounding $e^{-ax}\ge e^{-\pi}$ on $0<x\le\pi/|\xi|$ when $|\xi|\ge a$,
\begin{equation}
  I_a(\xi)\ \ge\ e^{-\pi}\frac{2\xi^2}{\pi^2}\int_0^{\pi/|\xi|}x^{1-Y}dx
  =\frac{2e^{-\pi}}{(2-Y)\pi^{Y}}\,|\xi|^{Y}=\kappa_Y|\xi|^Y\quad(|\xi|\ge a).
  \label{eq:LK-high}
\end{equation}
Both bounds hold for every $Y\in(0,2)$; the hypothesis $Y>1$ is needed only to
control the drift by the $H^{Y/2}$ norm in Step~2.

\medskip\noindent
\textbf{Step 4: the weight.}
The weight of \eqref{eq:HMG-norm} with $q=s_\mu=Y/2$ is
$w_\mu(\xi)=1+(M^2+\xi^2)^{Y/2}+(G^2+\xi^2)^{Y/2}$.  Since
$(a^2+\xi^2)^{Y/2}\le 2^{Y/2}\max\{a,|\xi|\}^{Y}$,
\begin{equation}
  w_\mu(\xi)\le 1+2^{1+Y/2}a_\ast^{Y}\quad(|\xi|\le a_\ast),\qquad
  w_\mu(\xi)\le 1+2^{1+Y/2}|\xi|^{Y}\quad(|\xi|\ge a_\ast).
  \label{eq:weight-bounds}
\end{equation}

\medskip\noindent
\textbf{Step 5: combination.}
For $|\xi|\le a_\ast$, \eqref{eq:uniform-low-freq} and \eqref{eq:weight-bounds} give
$\Real A_\mu(\xi)\ge \rf\,(1+2^{1+Y/2}a_\ast^{Y})^{-1}w_\mu(\xi)$.
For $|\xi|\ge a_\ast\ge\max\{G,M\}$, \eqref{eq:LK-high} applies to both $I_G$ and $I_M$, so
$\Real A_\mu(\xi)\ge \rf+2C\kappa_Y|\xi|^Y$, and with \eqref{eq:weight-bounds} and
$(\alpha+\beta)/(\gamma+\delta)\ge\min\{\alpha/\gamma,\beta/\delta\}$,
\[
  \frac{\Real A_\mu(\xi)}{w_\mu(\xi)}\ \ge\ \min\Bigl\{\rf,\ \frac{2C\kappa_Y}{2^{1+Y/2}}\Bigr\}
  \ \ge\ \min\Bigl\{\frac{\rf}{1+2^{1+Y/2}a_\ast^{Y}},\ \frac{C\kappa_Y}{2^{Y/2}}\Bigr\}.
\]
Hence $\Real A_\mu(\xi)\ge c_1\,w_\mu(\xi)$ for all $\xi$ with $c_1$ as in
\eqref{eq:c1-explicit}, and Plancherel's theorem gives
$\Real a_\R(u,u)\ge c_1\norm{u}_{H^{s_\mu}_{M,G}(\R)}^2$ directly in the weighted norm,
without the equivalence factor of Step~1.

\medskip\noindent
\textbf{Step 6: transfer to $\Omega$.}
For $v,w\in\Vmu$ the zero extensions satisfy
$\widetilde v,\widetilde w\in H^{s_\mu}(\R)$ with
$\norm{\widetilde v}_{L^2(\R)}=\norm{v}_{L^2(\Omega)}$, and by definition
$a_\mu(v,w)=a_\R(\widetilde v,\widetilde w)$ and
$\norm{v}_{\Vmu}=\norm{\widetilde v}_{H^{s_\mu}_{M,G}(\R)}$.  Both estimates
therefore transfer with identical constants, independent of $h$ and of
$\Omega$.  Finally, for real-valued $v$ the form $a_\mu(v,v)$ is real, so
$a_\mu(v,v)=\Real a_\mu(v,v)$.  \qed

The constants are explicit admissible choices rather than optimal ones.
On the box of equation~(\ref{eq:numerical-box}) of the main text, $c_1=5.11\times10^{-5}$ (attained at
$Y=Y_{\max}$, $\rf=r_{\min}$, $a_\ast=M_{\max}=9$) and $C_{\rm cont}\approx2.5\times10^3$.
The infimum of $\Real A_\mu(\xi)/w_\mu(\xi)$ over the box is in fact attained at $\xi=0$, where
\eqref{eq:uniform-low-freq} is an equality: its value $r_{\min}/(1+M_{\max}^{Y_{\max}}+G_{\max}^{Y_{\max}})
=1.40\times10^{-4}$ is confirmed numerically over $544$ box points and $\xi\in[0,10^6]$
(\texttt{results/coercivity\_c1.json}), so the explicit constant is within a factor three of the
optimal one.  The earlier version of this appendix obtained $c_1\approx1.6\times10^{-8}$ from a
high-frequency argument that required $Y_{\min}>1$ at every step; the L\'evy--Khintchine form
removes both the loss and that requirement from the coercivity part of the proof.

\section{Exact-assembly consequences and constrained dual norms}
\label{ssec:exact-assembly}
Under exact assembly, $z^\top\Sh z=a_\mu(z_h,z_h)$, so uniform coercivity
implies $\Sh\succ0$. The symmetric parts of $\Mh+\dt\Ah/2$ and
$\Mh+\dt\Ah$ are positive definite, giving nonsingularity.
\begin{lemma}[Positive-rate mass lower bound for exact assembly]
\label{lem:mass-lower}
Under exact assembly (Section~\ref{ssec:exact-assembly}), for
$\mu\in\Pset$ (so $\rf>0$),
\begin{equation}
  \Sh(\mu)\succeq \rf\,\Mh
  \label{eq:S-ge-rM}
\end{equation}
in the Loewner order, and consequently
\begin{equation}
  \bb^\top\Sh(\mu)^{-1}\bb\le \rf^{-1}\,\bb^\top\Mh^{-1}\bb
  \qquad\text{for every }\bb\in\R^{N_h}.
  \label{eq:mass-upper}
\end{equation}
\end{lemma}

\begin{proof}
By Lemma~\ref{lem:periodisation} and Plancherel's theorem,
\[
  \zz^\top\Sh\zz=a_\mu(z_h,z_h)=a_\R(\widetilde z_h,\widetilde z_h)
  =\frac1{2\pi}\int_\R\Real A_\mu(\xi)\,\bigl|\widehat{\widetilde z_h}(\xi)\bigr|^2\,d\xi .
\]
The exact symbol identity of \cite{DamircheliCGMY2026}, namely
\begin{equation}
\Real A_\mu(\xi)
=\rf+C\Gamma(-Y)\Bigl[\bigl(M^Y-\Real(M-i\xi)^Y\bigr)
+\bigl(G^Y-\Real(G+i\xi)^Y\bigr)\Bigr],
\label{eq:symbol-real-exact}
\end{equation}
together with $\Gamma(-Y)>0$ for $Y\in(1,2)$ and the monotonicity
$\Real(a\pm i\xi)^Y\le a^Y$, gives $\Real A_\mu(\xi)\ge \rf$ for all $\xi$.  Hence
\[
  \zz^\top\Sh\zz\ \ge\ \rf\norm{\widetilde z_h}^2_{L^2(\R)}
  =\rf\norm{z_h}^2_{L^2(\Omega)}=\rf\,\zz^\top\Mh\zz .
\]
Inversion of the Loewner
inequality between positive definite matrices reverses it, giving
\eqref{eq:mass-upper}.
\end{proof}
Exact assembly also implies the mesh-uniform bound on $\kappa_h$:
\[
 \|\Ah v\|_{\Sh^{-1}}
 =\sup_{w\ne0}\frac{a_\mu(v_h,w_h)}{\|w\|_{\Sh}}
 \le\frac{C_{\rm cont}}{c_1}\|v\|_{\Sh}.
\]
For any assembled matrix with SPD symmetric part, duality gives
$\|\Ah v\|_{\Sh^{-1}}\ge v^\top\Ah v/\|v\|_{\Sh}=\|v\|_{\Sh}$,
proving $\kappa_h\ge1$. This quantity is an operator norm, rather than the
condition number of $\Sh$.
\begin{proposition}[Conservatism of a full-space dual norm]
\label{prop:full-space-conservative}
Suppose $\bar S\succ0$, and the reference and candidate have identical
prescribed coefficients. Then $S_I\succ0$ and
\begin{equation}
 \rho_I^\top S_I^{-1}\rho_I
 \le (E\rho_I)^\top\bar S^{-1}(E\rho_I).
 \label{eq:full-space-conservative}
\end{equation}
Consequently a valid upper Riesz bracket for the zero-padded free residual
using $\bar S$ is also a conservative reliability bound for the constrained
problem, subject to the same initial-error and arithmetic assumptions.
\end{proposition}
\begin{proof}
The dual quadratic form equals
$\sup_v(2r^\top v-v^\top S v)$. Restricting the full-space supremum to
$v=Ew$ yields the free-space supremum, proving
\eqref{eq:full-space-conservative}. Positive definiteness restricts to every
principal subspace.
\end{proof}
The converse positivity implication fails because prescribed directions
can make $\bar S$ indefinite while $S_I$ is positive. The entire residual,
Riesz action and spectral bound must use the same constrained space.

\section{Discrete parameter stability}
\label{ssec:solution-operator}
\begin{proposition}[Discrete parameter stability of the solution map]
\label{prop:discrete-lipschitz}
Fix a mesh, a time grid, and a payoff $\psi$; the assembled symmetric
matrix $\Sh(\mu_1)$ is positive definite by Theorem~\ref{thm:spectral}. Let
$\mu_1,\mu_2\in\Pset$
with $\cc_h(\mu_1),\cc_h(\mu_2)$ the corresponding Crank--Nicolson
trajectories of \eqref{eq:cn-reference}.  Then, for every
$m=1,\dots,N_T$,
\begin{equation}
\begin{aligned}
  \norm{\cc_h^m(\mu_1)-\cc_h^m(\mu_2)}_{\Mh}^2
  \le{}& \norm{\cc_h^0(\mu_1)-\cc_h^0(\mu_2)}_{\Mh}^2\\
  &+\dt\sum_{n=0}^{m-1}
  \bigl\lVert\bm\varrho^{\,n+1/2}\bigr\rVert_{\Sh(\mu_1)^{-1}}^2 ,
\end{aligned}
  \label{eq:discrete-lipschitz}
\end{equation}
where the parameter-perturbation residual is
\begin{equation}
  \bm\varrho^{\,n+1/2}
  :=\bigl(\Ah(\mu_1)-\Ah(\mu_2)\bigr)\overline\cc_h^{\,n+1/2}(\mu_2)
   -\bigl(\overline\ff^{\,n+1/2}(\mu_1)-\overline\ff^{\,n+1/2}(\mu_2)\bigr).
  \label{eq:param-residual}
\end{equation}
Under exact assembly of the conforming variational matrices,
Theorem~\ref{thm:uniform-coercivity} additionally gives
\begin{equation}
  \norm{\bigl(\Ah(\mu_1)-\Ah(\mu_2)\bigr)\vv}_{\Sh(\mu_1)^{-1}}
  \le
  \frac{1}{c_1}\,
  \sup_{0\ne v_h,w_h\in\Vh}
  \frac{|a_{\mu_1}(v_h,w_h)-a_{\mu_2}(v_h,w_h)|}
       {\norm{v_h}_{\Vp{\mu_1}}\norm{w_h}_{\Vp{\mu_1}}}\,
  \norm{\vv}_{\Sh(\mu_1)},
  \label{eq:matrix-difference}
\end{equation}
so a bound on the difference of the two bilinear forms in a \emph{single}
norm---that of $\mu_1$---suffices, and no common energy space is needed.
\end{proposition}

\begin{proof}
The difference $\ee^n:=\cc_h^n(\mu_1)-\cc_h^n(\mu_2)$ satisfies the
Crank--Nicolson error equation \eqref{eq:discrete-error-equation} with
$\rr_\theta^{n+1/2}$ replaced by $\bm\varrho^{\,n+1/2}$ of
\eqref{eq:param-residual}, since $\cc_h(\mu_2)$ solves the $\mu_2$ scheme
exactly.  Apply Theorem~\ref{thm:fully-discrete-certificate} with
$\Sh=\Sh(\mu_1)$.  Estimate \eqref{eq:matrix-difference} follows from
the continuity argument above applied to the bilinear form
$a_{\mu_1}-a_{\mu_2}$.
\end{proof}
This fixed-grid statement compares both trajectories in the norm of the
first parameter. A continuum regularity result across varying $Y$ would
require additional control of the changing energy spaces.

\section{Other time schemes and exponential shifts}
\label{ssec:semidiscrete}\label{ssec:be-certificate}\label{ssec:shifted}
The following extensions follow the same energy argument as
Theorem~\ref{thm:fully-discrete-certificate}.
For a differentiable candidate and a semidiscrete reference, define
$\rho=\Mh\dot\cc_\theta+\Ah\cc_\theta-\ff$. Then
\begin{equation}
 \|\ee(t)\|_{\Mh}^2+\int_0^t\|\ee\|_{\Sh}^2
 \le\|\ee(0)\|_{\Mh}^2+\int_0^t\|\rho\|_{\Sh^{-1}}^2.
 \label{eq:semidiscrete-certificate}
\end{equation}
Multiply the error equation by $\ee^\top$, apply Young's inequality, and
integrate.
For backward Euler, the residual is
$\rho^{n+1}=\Mh(c_\theta^{n+1}-c_\theta^n)/\dt_n
+\Ah c_\theta^{n+1}-f^{n+1}$ and
\begin{equation}
 \|\ee^{n+1}\|_{\Mh}^2+\|\ee^{n+1}-\ee^n\|_{\Mh}^2
 +\dt_n\|\ee^{n+1}\|_{\Sh}^2
 \le\|\ee^n\|_{\Mh}^2+\dt_n\|\rho^{n+1}\|_{\Sh^{-1}}^2.
 \label{eq:be-certificate}
\end{equation}
This follows from
$2y^\top\Mh(y-z)=\|y\|_{\Mh}^2-\|z\|_{\Mh}^2+\|y-z\|_{\Mh}^2$.
A Rannacher start combines backward-Euler half steps with CN
\cite{Rannacher1984}; summing the corresponding bounds, with the initial
term retained, gives the certificate for the mixed trajectory.
Each step must use its own scheme's residual.

For an exact-assembly G\aa rding bound
$a_\mu(v,v)\ge c_1\|v\|_{V_\mu}^2-c_2\|v\|_{L^2}^2$, choose
$\sigma>c_2$ and set
$A_\sigma=\Ah+\sigma\Mh$, $S_\sigma=\Sh+\sigma\Mh$.
Then $S_\sigma\succeq(\sigma-c_2)\Mh$.
Advance $d_h$ by CN for the shifted operator with forcing
$g^n=e^{-\sigma\tau_n}f^n$, and reconstruct
$c_{h,\sigma}^n=e^{\sigma\tau_n}d_h^n$.
For $d_\theta^n=e^{-\sigma\tau_n}c_\theta^n$ and its shifted CN residual
$\rho_\sigma$, the same theorem yields
\[
 \max_m\|c_{h,\sigma}^m-c_\theta^m\|_{\Mh}
 \le e^{\sigma T}
 \left(\|d_h^0-d_\theta^0\|_{\Mh}^2+
 \dt\sum_n\|\rho_\sigma^{n+1/2}\|_{S_\sigma^{-1}}^2\right)^{1/2}.
\]
This certifies the reconstructed shifted target. Exponential transformation
and CN discretisation do not commute, so it is a different discrete
trajectory from unshifted CN. The numerical studies use unshifted CN with
positive rates.

\section{Derivative bounds and continuum error separation}
\label{ssec:derivatives}\label{ssec:total-error}
For an $H^2(\Omega)$ spline space (cardinal order $p\ge3$), define the matrices
$R_h^{(j)}=((\phi_{h,k'}^{(j)},\phi_{h,k}^{(j)}))$ and the constants
$\Lambda_j=\lambda_{\max}(R_h^{(j)},\Mh)$ for $j=1,2$.
Since the common lifting cancels,
\[
 \|\partial_x^j(V_h^m-V_{\theta,h}^m)\|_{L^2}
 \le\sqrt{\Lambda_j}\eta_{\rm aug},
 \qquad \Lambda_j=O(h^{-2j}).
\]
Indeed, $\ee^\top R_h^{(j)}\ee\le\Lambda_j\ee^\top\Mh\ee$.
With the Greeks $\Delta=e^{-x}V_x$ and $\Gamma_{\rm opt}=e^{-2x}(V_{xx}-V_x)$,
\[
 \|\Delta_h-\Delta_\theta\|_{L^2}
 \le e^{-x_{\min}}\sqrt{\Lambda_1}\eta_{\rm aug},\qquad
 \|\Gamma_{{\rm opt},h}-\Gamma_{{\rm opt},\theta}\|_{L^2}
 \le e^{-2x_{\min}}(\sqrt{\Lambda_2}+\sqrt{\Lambda_1})\eta_{\rm aug}.
\]
These global factors can be conservative on a wide log-price domain and
grow under refinement. A rigorous numerical supremum of the point factor
in \eqref{eq:gamma-infty} requires validated cellwise extrema.

Under exact assembly, the regularity and buffer assumptions of
\cite{DamircheliCGMY2026}, and an identical prescribed smooth lifting,
\begin{equation}
 \max_n\|V(\tau_n)-V_h^n\|_{L^2(\Omega)}
 \le C_\mu(h^{p-Y/2}+\dt^2+\varepsilon_b).
 \label{eq:deterministic-error}
\end{equation}
The lifting discrepancies vanish because it is prescribed identically.
Combining this estimate with \eqref{eq:total-error-split} gives
\[
 \|V(\tau_m)-V_{\theta,h}^m\|_{L^2}
 \le C_\mu(h^{p-Y/2}+\dt^2+\varepsilon_b)+\eta_{\rm aug}.
\]
The deterministic constants and buffer error are not evaluated as
computable bounds here; the first term is assessed instead by the price comparison of
Section~\ref{ssec:calibration}, and the certificate supplies the second.

\section{Residual certificate for the American put}
\label{ssec:vi}
For a discrete obstacle problem in coefficient space (B-splines are
nonnegative, so $c\ge\psi_h$ coefficientwise implies $U_h\ge\psi_h$), one
backward-Euler step reads: find $u^{n+1}\ge\psi_h$ with
$(B_+u^{n+1}-b^n)^\top(v-u^{n+1})\ge0$ for all $v\ge\psi_h$, where
$B_+=\Mh+\dt\Ah$ and $b^n=\Mh u^n+\dt f^{n+1}$.
\begin{theorem}[Residual certificate for the discrete variational inequality]
\label{thm:vi}
Let $\{u_\theta^n\}$ be any candidate with $u_\theta^n\ge\psi_h$ for all $n$,
let $R^n:=B_+u_\theta^{n+1}-\Mh u_\theta^n-\dt f^{n+1}$, and split
$R^n=R^n_+-R^n_-$ into positive and negative parts.  Then for every $m$,
\[
 \norm{e^m}_{\Mh}^2+\dt\sum_{n<m}\norm{e^{n+1}}_{\Sh}^2
 \ \le\ \norm{e^0}_{\Mh}^2
 +\sum_{n<m}\Bigl[\,2\,(R^n_+)^\top(u_\theta^{n+1}-\psi_h)
 +\dt^{-1}\norm{R^n_-}_{\Sh^{-1}}^2\Bigr].
\]
\end{theorem}
\begin{proof}
Test the inequality with $v=u_\theta^{n+1}$:
\[
(B_+e^{n+1}-\Mh e^n,\,e^{n+1})\le-(R^n,e^{n+1})=(R^n,u_\theta^{n+1}-u^{n+1}).
\]
Since $u^{n+1}\ge\psi_h$ and $R_+\ge0$, $(R^n_+,u_\theta^{n+1}-u^{n+1})\le(R^n_+,u_\theta^{n+1}-\psi_h)$,
and $(R^n_-,e^{n+1})\le\norm{R^n_-}_{\Sh^{-1}}\norm{e^{n+1}}_{\Sh}$.  The left side equals
\[
\tfrac12\bigl(\norm{e^{n+1}}_{\Mh}^2-\norm{e^n}_{\Mh}^2+\norm{e^{n+1}-e^n}_{\Mh}^2\bigr)+\dt\norm{e^{n+1}}_{\Sh}^2;
\]
Young's inequality and summation finish the proof.
\end{proof}
Everything in Sections~\ref{sec:galerkin}--\ref{sec:affine} (exact assembly, $\smin$, the banded norm,
the affine $A_Q$) applies verbatim; only the right-hand side gains the
complementarity term, which vanishes where the candidate is either off the
obstacle with nonpositive residual or on it.  The Crank--Nicolson version
requires testing the inequality at the midpoint, which is not admissible in
general; this is the one open technical point and the reason to use backward
Euler (or Rannacher-type restarts) for an American reference.

\section{Preconditioning and correction}
\label{ssec:correction-extension}
If $P_h\succ0$ and $\Sh\succeq\gamma_P P_h$ with $\gamma_P>0$, then
$r^\top\Sh^{-1}r\le\gamma_P^{-1}r^\top P_h^{-1}r$.
Alternatively $P_h$ can accelerate the Riesz iteration while the bracket
uses its independent bound on $\lambda_{\min}(\Sh)$.
\begin{proposition}[Contraction of a preconditioned correction]
\label{prop:correction-contraction}
Fix a time step and write the step equation as $B_+\zz=\bb$.  Let $P\succ0$ be
a preconditioner and consider the Richardson iteration
$\zz^{(i+1)}=\zz^{(i)}+P^{-1}(\bb-B_+\zz^{(i)})$ with residuals
$\sres^{(i)}=\bb-B_+\zz^{(i)}$.  Then
\begin{equation}
  \norm{\sres^{(i+1)}}_{P^{-1}}\le\rho_P\,\norm{\sres^{(i)}}_{P^{-1}},
  \qquad
  \rho_P:=\bigl\lVert I-P^{-1/2}B_+P^{-1/2}\bigr\rVert_2 .
  \label{eq:contraction}
\end{equation}
If $P=\Mh+\tfrac{\dt}{2}\Sh$, the symmetric part of $B_+$, then
\[
I-P^{-1/2}B_+P^{-1/2}=-\tfrac{\dt}{2}P^{-1/2}N_hP^{-1/2},\qquad N_h=\tfrac12(\Ah-\Ah^\top),
\]
with $N_h$ skew symmetric, so $\rho_P$ is determined entirely by the nonsymmetric part of the
CGMY operator and satisfies $\rho_P\le\tfrac{\dt}{2}\lVert P^{-1/2}N_hP^{-1/2}\rVert_2$.
\end{proposition}

\begin{proof}
$\sres^{(i+1)}=(I-B_+P^{-1})\sres^{(i)}$, so
$P^{-1/2}\sres^{(i+1)}=(I-P^{-1/2}B_+P^{-1/2})P^{-1/2}\sres^{(i)}$ and
\eqref{eq:contraction} follows by taking Euclidean norms.  For the second
claim, $B_+=P+\tfrac{\dt}{2}N_h$ by \eqref{eq:Sh} and the definition of $N_h$.
\end{proof}
The proposition proves contraction only when $\rho_P<1$; it is a
Richardson estimate and makes no contraction or timing assertion about the
GMRES correction used in the numerical experiments.
\FloatBarrier
\section{Verification of the assembly, the bounds and the solvers}
\label{sm:verification}
The studies of this section use the design of Section~\ref{ssec:design} of the main text (60 Crank--Nicolson
steps, several $N$); they verify the assembly of Section~\ref{ssec:assembly}, the bounds of
Theorems~\ref{thm:spectral} and~\ref{thm:comparison} and the solvers of Section~\ref{ssec:comparison}.
\subsection{Assembly verification}
Table~\ref{tab:assembly-v9} compares the assembled free matrix with the exact Galerkin matrix,
computed independently at 40 digits (alias series summed in closed form through Hurwitz zeta
tails, Toeplitz entries by a high-precision trapezoid rule).  The truncation error decreases
like $K^{-(2p-Y-1)}$ as predicted by Lemma~\ref{lem:tail}, and the bound $\varepsilon_K$ is
within a factor $1.1$--$3.7$ of the observed error at $K=100$ (at most $4.4$ at $K=20$ and $8.6$ at $K=5$).
\begin{table}[!htbp]
\centering\small
\caption{Assembly verification against the exact Galerkin matrix (40-digit real-line fixture): relative operator-norm error of the free matrix $A_I$, median / max over the ten cost-study queries ($N\le256$; three at $N=512$). The principal-period rule of the earlier version is the alias sum with $K=0$; $\varepsilon_K/\|A_I\|_2$ is the bound of Lemma~\ref{lem:tail}.}
\label{tab:assembly-v9}
\resizebox{\textwidth}{!}{%
\begin{tabular}{rrrrrrr}
\toprule
$N$ & queries & principal period & $K=5$ & $K=20$ & $K=100$ & $\varepsilon_{100}/\|A_I\|_2$ (max)\\
\midrule
64 & 10 & 0.381 / 0.410 & $3.0\times10^{-5}$ / $5.8\times10^{-5}$ & $2.8\times10^{-7}$ / $7.7\times10^{-7}$ & $9.7\times10^{-10}$ / $4.2\times10^{-9}$ & $4.8\times10^{-9}$\\
128 & 10 & 0.370 / 0.408 & $2.5\times10^{-5}$ / $5.4\times10^{-5}$ & $2.3\times10^{-7}$ / $7.0\times10^{-7}$ & $8.1\times10^{-10}$ / $3.8\times10^{-9}$ & $4.4\times10^{-9}$\\
256 & 10 & 0.358 / 0.404 & $2.3\times10^{-5}$ / $5.1\times10^{-5}$ & $2.0\times10^{-7}$ / $6.6\times10^{-7}$ & $7.1\times10^{-10}$ / $3.6\times10^{-9}$ & $4.1\times10^{-9}$\\
512 & 3 & 0.333 / 0.402 & $1.6\times10^{-5}$ / $4.9\times10^{-5}$ & $1.2\times10^{-7}$ / $6.4\times10^{-7}$ & $3.6\times10^{-10}$ / $3.5\times10^{-9}$ & $4.0\times10^{-9}$\\
\bottomrule
\end{tabular}}
\end{table}

\subsection{Spectral bound and comparison norm}
\begin{figure}[!htbp]
\centering
\includegraphics[width=0.94\linewidth]{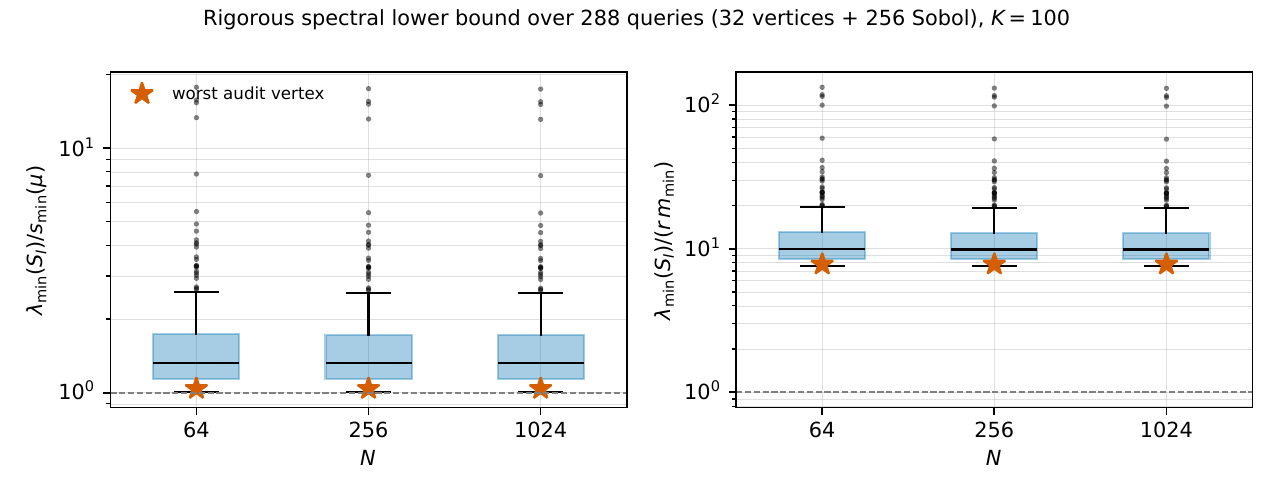}
\caption{Ratios $\lambda_{\min}(S_I)/\smin(\mu)$ (left) and $\lambda_{\min}(S_I)/(\rf m_{\min})$ (right)
over the 288 verification queries at $N=64,256,1024$; the star marks the worst vertex of the
earlier positivity audit.  Both ratios are $\ge1$ at every query, as Theorem~\ref{thm:spectral} asserts.}
\label{fig:spectral-bound}
\end{figure}
\begin{figure}[!htbp]
\centering
\includegraphics[width=0.94\linewidth]{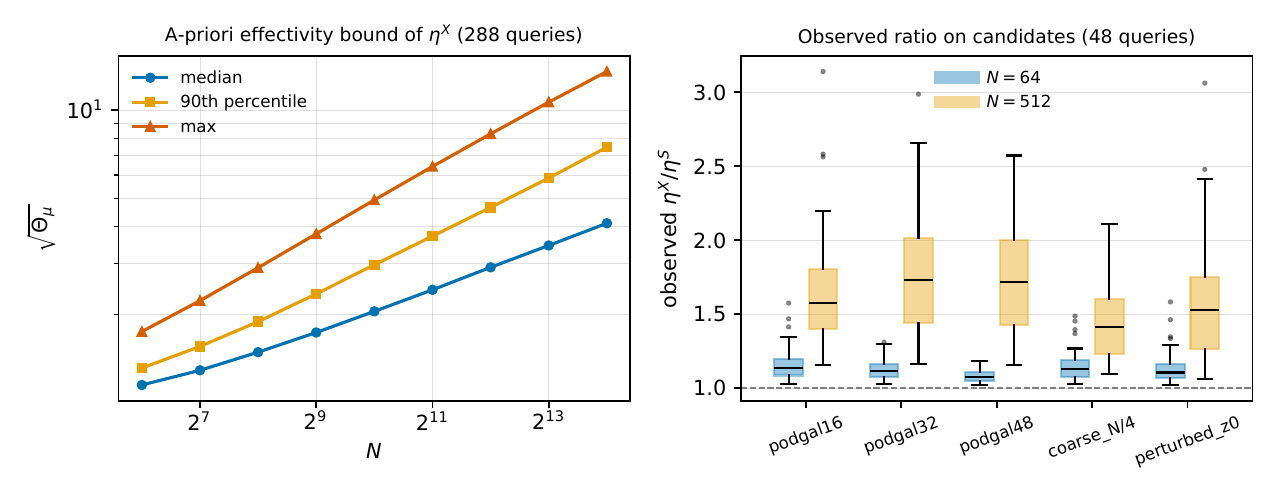}
\caption{Left: the a-priori effectivity bound $\sqrt{\Thmu}$ of the banded certificate as a function of
$N$ (symbol only, 288 queries).  Right: observed ratio $\eta^X/\eta^S$ for POD--Galerkin candidates
(ranks 16, 32, 48), the prolonged coarse solution and exact dynamics from a perturbed initial vector,
at $N=64$ and $N=512$.}
\label{fig:theta}
\end{figure}
The bound of Theorem~\ref{thm:spectral} and the two-sided comparison \eqref{eq:comparison} hold at
all $289\times3$ tested (query, $N$) pairs (the 288 verification queries and the worst audit vertex); the generalised eigenvalues of $(S_I,X_I)$ lie in
$[1,\Thmu]$ in every case.  The a-priori factor $\sqrt{\Thmu}$ grows like
$(\min(G,M)h)^{(Y-2)/2}$ (median $1.15$, $1.74$, $2.05$, $2.91$, $4.11$ at $N=64$, $512$, $1024$,
$4096$, $16384$; maxima $1.75$, $3.78$, $4.93$, $8.30$, $13.6$), because the tempering scale moves
below the resolution of a degree-$(p-1)$ trigonometric minorant.  The observed ratio
$\eta^X/\eta^S$ on candidates is smaller: median $1.08$--$1.13$ at $N=64$ and $1.4$--$1.7$ at
$N=512$ (maximum $3.1$), so at the production resolution the banded certificate costs a factor
below two in effectivity; the production candidates of Section~\ref{ssec:e1} confirm the factor
$\eta^X/\eta^S\approx2$ at $N=1024$.
\subsection{Riesz solves and reference trajectories}
\begin{table}[!htbp]
\centering\small
\caption{Iteration counts (median / max over 20 Sobol queries) to relative residual $10^{-12}$: conjugate gradients for the Riesz system $S_Iv=\rho$ with and without the circulant preconditioner built from the symbol samples, and circulant-preconditioned GMRES per Crank--Nicolson step for $B_+=M_I+\tfrac{\Delta\tau}{2}A_I$ (60 steps, warm start, mean over steps / max).}
\label{tab:pcg-v9}
\begin{tabular}{rrrr}
\toprule
$N$ & plain CG & circulant PCG & circulant GMRES ($B_+$)\\
\midrule
64 & 42 / 45 & 8 / 8 & 5.0 / 5\\
128 & 74 / 82 & 9 / 9 & 5.0 / 6\\
256 & 130 / 154 & 9 / 10 & 6.0 / 7\\
512 & 224 / 292 & 10 / 10 & 7.0 / 8\\
1024 & 384 / 548 & 10 / 11 & 7.1 / 12\\
2048 & 648 / 1025 & 11 / 11 & 8.0 / 13\\
4096 & 1080 / 1915 & 11 / 12 & 9.0 / 15\\
\bottomrule
\end{tabular}
\end{table}

Dense-LU and circulant-GMRES references agree to $10^{-12}$ in the $M_I$ norm for
$N\le1024$; the bracket width at convergence is below $10^{-20}$ relative for $N\le1024$ and
below $10^{-18}$ at $N=4096$, where it is set by the roundoff of the recomputed residual.
\subsection{Initial-error tests and \texorpdfstring{$Y<1$}{Y<1}}
On the 48 exact-dynamics cases with perturbed initial vectors ($\|\delta\|_{M_I}=10^{-4}$), both
$\eta^X_{\rm aug}$ and $\eta^S_{\rm aug}$ hold with effectivity $1.00\pm10^{-6}$, and the
residual-only estimators underestimate the error in every case.  For $Y\in\{0.5,0.8,0.95\}$
(same $C,G,M,\rf$ as the reference query) the bound $\alpha_{\rm lb}=0.0500=\rf$ is positive,
$\lambda_{\min}(S_I)\in[0.052,0.055]$, $\sqrt{\Thmu}\le1.5$, the preconditioned GMRES solve of
$B_+$ needs at most five iterations per step, and the certificate holds with effectivity one on
perturbed exact dynamics; the discrete framework therefore covers the finite-variation regime
without change.
\section{The neural-operator candidate of the earlier version}
\label{sm:deeponet}
The neural-operator candidate is a DeepONet with payoff and parameter branches and a time
trunk, latent dimension 64 and 80,002 trainable parameters, followed by a coefficient decoder
that supplies the spline-valued output. The experiments fix $\mathcal K=100$;
functional-input notation accommodates a payoff encoder but the benchmark
tests a parameter map for this single strike.  It was trained on the principal-period discretisation of the earlier version of this work, and
we did not retrain it on the exact-assembly references.  Certified as a candidate of the exact
system ($N=64$, $N_T=60$, 120 benchmark queries), it is rejected at every tolerance of the sweep
$10^{-3}\times\{0.1,\dots,10\}$: the median true error is $0.78$ (relative $3.1\times10^{-3}$), the
median exact bound $3.4$, with no false accept.  The certificate thus measures, correctly, the
distance between two discrete problems; the neural candidate of the production study is the
POD-NN below.
\section{Neural reduced candidates: results and training details}
\label{sm:training}
Table~\ref{tab:podnn} reports every POD-NN checkpoint of Section~\ref{ssec:e1} of the main text on
the 30 validation queries, and Table~\ref{tab:training} lists the optimiser, schedule, parameter
counts and wall times of the runs.
\begin{table}[!htbp]
\centering\small
\caption{Neural reduced candidates (POD-NN, $r=48$) on the 30 validation queries at the production configuration: median true error and median augmented bounds of the three certificates, with the effectivity of the exact one. Standard output: $r\times N_T$ coefficients; space--time output: $r\times r_t$ coefficients of a temporal POD ($r_t=23$). Checkpoints selected by the validation median of the reduced certificate. The projection of the reference trajectories onto the same basis and the POD--Galerkin solution are the attainable floor. No candidate passes at $\varepsilon_N=1.18\times10^{-3}$; every POD-NN query is settled by the cheap rejection.}
\label{tab:podnn}
\resizebox{\textwidth}{!}{%
\begin{tabular}{lrrrrr}
\toprule
candidate & error & $\eta^S_{\rm aug}$ & $\eta^X_{\rm aug}$ & $\eta^{\bar X}_{\rm aug}$ & eff.\ $\eta^S$\\
\midrule
800 traj., $\lambda_R=0$ & 0.183 & 33.7 & 42.0 & 103 & 178.5\\
800 traj., $\lambda_R=1$ & 1.76 & 5.20 & 8.19 & 18.9 & 3.0\\
4000 traj., $\lambda_R=0$ & 0.0816 & 12.2 & 16.4 & 37.4 & 148.0\\
4000 traj., $\lambda_R=10^{-2}$ & 0.448 & 1.27 & 2.11 & 4.90 & 2.9\\
4000 traj., $\lambda_R=10^{-1}$ & 0.485 & 1.28 & 2.15 & 4.95 & 2.8\\
4000 traj., $\lambda_R=1$, seed 0 & 0.458 & 1.26 & 2.07 & 4.94 & 2.9\\
\quad seed 1 & 0.473 & 1.24 & 2.04 & 4.80 & 2.8\\
\quad seed 2 & 0.464 & 1.27 & 2.08 & 4.95 & 2.8\\
space--time, 4000 traj., $\lambda_R=1$, seed 0 & 0.0159 & 0.0582 & 0.0891 & 0.190 & 3.7\\
\quad seed 1 & 0.0161 & 0.0559 & 0.0841 & 0.195 & 3.6\\
\quad seed 2 & 0.0142 & 0.0545 & 0.0815 & 0.176 & 3.8\\
\midrule
projection of the references & $3.3\times10^{-4}$ & $3.2\times10^{-3}$ & $6.2\times10^{-3}$ & 0.0127 & 10.2\\
POD--Galerkin, $r=48$ & $3.5\times10^{-4}$ & $3.8\times10^{-3}$ & $8.2\times10^{-3}$ & 0.0154 & 11.0\\
\bottomrule
\end{tabular}}
\end{table}

\begin{table}[!htbp]
\centering\small
\caption{POD-NN training details (production configuration). MLP $6\to256\to256\to256\to$ output with SiLU, double precision, inputs $(\mu,T)$ standardised, $\alpha^0$ fixed to the projection of $z_h^0$; loss = relative least squares in the reduced coordinates $+\lambda_R\,\eta^{\bar X}$ (reduced online certificate, differentiable); Adam, learning rate $10^{-3}$ with cosine decay to $10^{-5}$, 2000 epochs, batch 64, validation every 10 epochs on the 30-query set, checkpoint = epoch with the smallest median reduced bound. Intel Core Ultra 9 285K, 3 threads per run.}
\label{tab:training}
\resizebox{\textwidth}{!}{%
\begin{tabular}{lrrrr}
\toprule
variant & output & parameters & runs & wall time per run (h)\\
\midrule
standard, 800 trajectories & $r\times N_T=11\,520$ & $3.09$~M & $\lambda_R\in\{0,10^{-2},10^{-1},1\}$, seed 0 & 0.9--2.7\\
standard, 4000 trajectories & $r\times N_T=11\,520$ & $3.09$~M & $\lambda_R\in\{0,10^{-2},10^{-1},1\}$, seed 0; $\lambda_R=1$, seeds 1--2 & 3.5--7.8\\
space--time, 4000 trajectories & $r\times r_t=1\,104$ & $0.42$~M & $\lambda_R=1$, seeds 0--2 & 5.9--7.1\\
\bottomrule
\end{tabular}}
\end{table}

\section{Measured constants of the cost model}
\label{sm:cost}
Table~\ref{tab:cost} gives the measured constants of the cost model of Section~\ref{ssec:cost} of the
main text, from the refinement study of Section~\ref{ssec:benchmark}.
\begin{table}[!htbp]
\centering\small
\caption{Cost model of Section~\ref{ssec:cost} with measured constants (one thread, idle machine; medians over five benchmark queries; $m=2N_h$, $i_n$ and $p_n$ the GMRES and PCG iteration counts per step). Times in ms; $c$ is the fitted time per operation of the stated count with its maximum relative deviation; where the per-step interpreter overhead matters the two-constant fit $t=c'\,\mathrm{ops}+d\,\mathrm{calls}$ is also given. The spectral-data row is the fine verification grid; the online path uses the fast mode ($2.5$~ms, Section~\ref{ssec:benchmark}). The reduced rows do not depend on $N$ by construction.}
\label{tab:cost}
\resizebox{\textwidth}{!}{%
\begin{tabular}{lrrrrrrr}
\toprule
component & op count & $N=256$ & $N=512$ & $N=1024$ & $N=2048$ & $N=4096$ & $c$ (fit res.)\\
\midrule
assembly & $L(2K+1)$ & 21.63 & 57.02 & 106.37 & 209.78 & 408.07 & $1.25\times10^{-7}$ (16\%)\\
spectral data (fine mode; online fast mode $2.5$~ms) & $n_g(2K+1)$ & 477.85 & 478.78 & 476.70 & 476.52 & 476.85 & $1.23\times10^{-7}$ (1\%)\\
direct solve & $(\sum_n i_n+N_T)\,m\log m$ & 76.64 & 78.41 & 148.17 & 246.94 & 591.73 & $4.52\times10^{-9}$ (156\%); $c'=2.2\times10^{-9}$, $d=3.7\times10^{-5}$ (11\%)\\
residual & $N_T\,m\log m$ & 1.11 & 2.47 & 6.18 & 15.27 & 53.96 & $1.25\times10^{-9}$ (68\%)\\
banded $\eta^X$ & $N_T N_h$ & 0.80 & 1.61 & 3.33 & 6.76 & 13.38 & $1.35\times10^{-8}$ (6\%); $c'=1.3\times10^{-8}$, $d=0$ (6\%)\\
exact $\eta^S$ & $\sum_n p_n\,m\log m$ & 61.49 & 60.83 & 141.40 & 237.32 & 678.46 & $2.65\times10^{-9}$ (115\%); $c'=1.7\times10^{-9}$, $d=1.2\times10^{-5}$ (21\%)\\
\midrule
\multicolumn{8}{l}{reduced model at $N=1024$, $Q=36$ (independent of $N$ by construction)}\\
component & op count & $r=16$ & $r=32$ & $r=48$ & $r=64$ & $r=80$ & $c$ (fit res.)\\
reduced solve & $(Q{+}2)r^2+r^3+3N_Tr^2$ & 0.78 & 0.86 & 0.97 & 1.11 & 1.29 & $5.37\times10^{-10}$ (638\%); $c'=8.8\times10^{-11}$, $d=3.3\times10^{-6}$ (4\%)\\
reduced certificate & $\rho(r{+}1)(Q{+}2{+}2N_T)$ & 0.27 & 0.61 & 1.10 & 1.80 & 2.74 & $8.73\times10^{-11}$ (36\%); $c'=7.6\times10^{-11}$, $d=3.9\times10^{-7}$ (4\%)\\
$\ell(\nu)+\varepsilon_Q$ & $161(2K+1)$ & 2.73 & 2.73 & 2.73 & 2.73 & 2.73 & $8.56\times10^{-8}$ (3\%)\\
$\rho$ (retained rows) & & 263 & 390 & 524 & 668 & 819 & \\
\bottomrule
\end{tabular}}
\end{table}